\documentclass[12pt]{article}
\usepackage{amsmath,amssymb}
\usepackage{bigdelim,multirow}
\usepackage[all]{xy}
\usepackage{color}
\usepackage{bm}
\usepackage{indentfirst}
\usepackage{amsthm}
\usepackage{enumerate}
\usepackage{framed}
\usepackage{array,arydshln}
\allowdisplaybreaks
\newtheoremstyle{mystyle}
    {}
    {}
    {\normalfont}
    {}
    {\bfseries}
    {}
    { }
    {}
\theoremstyle{mystyle}

\newtheorem{thm}{Theorem}
\newtheorem{lem}{Lemma}
\newtheorem*{prf}{Proof.}
\newtheorem{ex}{Example}

\newtheorem{defi}{Definition}

\makeatletter
\newif\if@borderstar

\def\bordermatrix{\@ifnextchar*{%
 \@borderstartrue\@bordermatrix@i}{\@borderstarfalse\@bordermatrix@i*}%
}
\def\@bordermatrix@i*{\@ifnextchar[{\@bordermatrix@ii}{\@bordermatrix@ii[()]}}
\def\@bordermatrix@ii[#1]#2{%
\begingroup
 \m@th\@tempdima8.75\p@\setbox\z@\vbox{%
 \def\cr{\crcr\noalign{\kern 2\p@\global\let\cr\endline }}%
 \ialign {$##$\hfil\kern 2\p@\kern\@tempdima & \thinspace %
  \hfil $##$\hfil && \quad\hfil $##$\hfil\crcr\omit\strut %
  \hfil\crcr\noalign{\kern -\baselineskip}#2\crcr\omit %
  \strut\cr}}%
 \setbox\tw@\vbox{\unvcopy\z@\global\setbox\@ne\lastbox}%
 \setbox\tw@\hbox{\unhbox\@ne\unskip\global\setbox\@ne\lastbox}%
 \setbox\tw@\hbox{%
  $\kern\wd\@ne\kern -\@tempdima\left\@firstoftwo#1%
  \if@borderstar\kern 2pt\else\kern -\wd\@ne\fi%
 \global\setbox\@ne\vbox{\box\@ne\if@borderstar\else\kern 2\p@\fi}%
 \vcenter{\if@borderstar\else\kern -\ht\@ne\fi%
  \unvbox\z@\kern -\if@borderstar2\fi\baselineskip}%
 \if@borderstar\kern-2\@tempdima\kern2\p@\else\,\fi\right\@secondoftwo#1 $%
 }\null \;\vbox{\kern\ht\@ne\box\tw@}%
\endgroup
}
\makeatother
\begin{document}

\title{Toric Construction and Chow Ring of Moduli Space of Quasi Maps from $\mathbb{P}^1$ with Two Marked Points to $\mathbb{P}^1 \times \mathbb{P}^1$}
\author{Kohki Matsuzaka\\
\it Division of Mathematics, Graduate School of Science \\
\it Hokkaido University \\
\it  Kita-ku, Sapporo, 060-0810, Japan \\
\it e-mail address: kohki@math.sci.hokudai.ac.jp
}
\date{}
\maketitle
\begin{abstract}
In this paper, we present explicit toric construction of moduli space of quasi maps from $\mathbb{P}^1$ with two marked points to $\mathbb{P}^1 \times \mathbb{P}^1$, which was first proposed by Jinzenji and prove that it is a compact orbifold. We also determine
its Chow ring and compute its Poincar\'e polynomial for some lower degree cases.  

\end{abstract}

\

\tableofcontents

\section{Introduction}

\subsection{The Moduli Space of Quasi Maps from $\mathbb{P}^1$ with Two Marked Points to $\mathbb{P}^1 \times \mathbb{P}^1$ and Its Compactification}

\textit{Mirror symmetry conjecture} is based upon existence of a pair of two Calabi--Yau varieties $(X, X^{\circ})$ such that  equivalence between A--model of $X$ (which depends on K\"{a}hler structure of $X$) and B--model of $X^{\circ}$ (which depends on complex structure of $X^{\circ}$) holds (see \cite{CK}). In the case of classical mirror symmetry,  it follows that A--model correlation function of $X$ and B--model correlation function of $X^{\circ}$ coincide (in case of compact Calabi-Yau 3-folds, these correlation functions are called Yukawa coupling). As the work \cite{CdOGP} by Candelas, de la Ossa, Green, and Parkes suggests, generating functions of Gromov--Witten invariants of the Calabi-Yau 3-fold $X$ (A--model Yukawa coupling) can be computed by using of B--model Yukawa coupling of $X^{\circ}$ obtained from analysis of Picard-Fuchs differential equation for period integral of the holomorphic 3-form of 
$X^{\circ}$. More precisely, the B-model Yukawa coupling of $X^{\circ}$ is translated into the A-model Yukawa coupling of $X$ via 
``the mirror map'', which is the coordinate transformation between the deformation parameter of K\"{a}hler structure of $X$ and 
the one of complex structure of $X^{\circ}$, obtained from solutions of the Picard-Fuchs equation.   
Then a question naturally arises. ``Can we interpret information of the solution of the Picard-Fuchs equation coming from 
the B--model of $X^{\circ}$ from the point of view of the A--model of $X$?''  The correlation functions of the A--model, 
Gromov-Witten invariantsm of $X$, are defined as intersection numbers of moduli space $\overline{M}_{0,3}(V,{\bf d})$
of stable maps of multi degree ${\bf d}$ from genus $0$ stable curve with three marked points to a toric variety $V$, 
which is an ambient space of $X$. At first sight, this moduli space seems to be far from the informations supplied by the B--model. 
But as the analysis of Gauged Linear Sigma model by Morrison and Plesser \cite{MP} tells us, another way of compactification of moduli space of holomorphic maps from $\mathbb{P}^{1}$ to $V$, i.e., so called toric compactification, is expected to supply information of the B--model.   

Following this expectation, Jinzenji conjectured existence of the moduli space $\widetilde{Mp}_{0,2} (V,\bf{d})$, which is a moduli space 
of quasimaps from $\mathbb{P}^{1}$ with two marked points to a toric variety $V$ of multi degree ${\bf d}$, compactified by quasimaps from nodal genus 0 curves \cite{J1}. The word ``quasi'' means that the map has several points on 
the nodal genus $0$ curve whose images are undefined. 
He also conjectured that generating function of intersection numbers constructed as analogue of Gromov-Witten invariants 
of $X$ reproduce the B--model information of $X^{\circ}$ including the mirror map, and confirmed the conjecture in 
several examples.
From these results, it is expected that $\widetilde{Mp}_{0,2} (V,\bf{d})$ is defined as a toric variety. In \cite{J1}, Jinzenji explicitly constructed homogeneous coordinates and weight matrix of $\mathbb{C}^{\times}$ actions on these coordinates in the cases that $V=\mathbb{P}^{N-1}
, \mathbb{P}^{1}\times \mathbb{P}^{1}$. But concrete toric constructions were not given even for these examples. 
Later, Saito constructed the toric data of $\widetilde{Mp}_{0,2} (\mathbb{P}^{N-1},d)$ \cite{S} and Jinzenji and Saito gave explicit 
toric construction of $\widetilde{Mp}_{0,2} (\mathbb{P}(3,1,1,1),d)$ \cite{JS}. On the other hand, explicit toric construction of $\widetilde{Mp}_{0,2} (\mathbb{P}^1 \times \mathbb{P}^1,\bf{d})$ has not been given.          

In this paper, we construct an explicit toric data of the moduli space $\widetilde{Mp}_{0,2} (\mathbb{P}^1 \times \mathbb{P}^1,\bf{d})$ 
and prove that the fan used in the construction is complete and simplicial, i.e., $\widetilde{Mp}_{0,2} (\mathbb{P}^1 \times \mathbb{P}^1,\bf{d})$  is a compact orbifold.  

First, we introduce the definition of $\widetilde{Mp}_{0,2} (\mathbb{P}^1 \times \mathbb{P}^1,\bf{d})$, where 
$\bf{d}$$=(d_1 ,d_2)$ \ $(d_1 \ge d_2 >0)$, which was first proposed in \cite{J1}.

\begin{defi}
\textit{
\begin{align}
& \widetilde{Mp}_{0,2} (\mathbb{P}^1 \times \mathbb{P}^1 , (d_1 ,d_2)) \notag \\
&:= \{ (\bm{a}_0 ,\dots ,\bm{a}_{d_1} ,\bm{b}_0 ,\dots ,\bm{b}_{d_2} , \widetilde{\bm{u}_0} ,\bm{u}_1 ,\dots ,\bm{u}_{d_1 -1} ,\widetilde{\bm{u}_{d_1}} ) \in \mathbb{C}^{(d_1 +3)(d_2 +3)-6} \  \big| \notag \\
&\  \bm{a}_0 ,\dots ,\bm{a}_{d_1},\bm{b}_0 ,\dots ,\bm{b}_{d_2} \in \mathbb{C}^{2}, \widetilde{\bm{u}_0} ,\widetilde{\bm{u}_{d_1}} \in \mathbb{C}^{d_2} , \bm{u}_1 ,\dots ,\bm{u}_{d_1 -1} \in  \mathbb{C}^{d_2 +1}, \notag \\
&\ \bm{a}_0 ,\bm{a}_{d_1},\bm{b}_0 ,\bm{b}_{d_2} \neq \bm{0}, \notag \\
&\ (\bm{a}_i , u_{(i,0)} \cdot \dots \cdot u_{(i,d_2)}) \neq \bm{0} \ (1 \le i \le d_1 -1), \notag \\  
&\ (\bm{b}_j , u_{(0,j)} \cdot \dots \cdot u_{(d_1 ,j)}) \neq \bm{0} \ (1 \le j \le d_2 -1), \notag \\
&\ (u_{(i_1 ,i_2)} , u_{(j_1 ,j_2)}) \neq (0,0) \notag \\
&\ \text{if} \ 0 \le i_1 < j_1 \le d_1 \ \text{and} \ 0 \le j_2 < i_2 \le d_2 \notag \\
&\ \} / (\mathbb{C}^{\times})^{(d_1 +1)(d_2 +1)+1} , \label{1.1.1}
\end{align}
where $\bm{a}_i =(a_i ^1 ,a_i ^2) \ (i=0, \dots ,d_1)$, $\bm{b}_j =(b_j ^1 ,b_j ^2) \ (j=0, \dots ,d_2)$, $\widetilde{\bm{u}_0}=(u_{(0,1)} ,u_{(0,2)}, \cdots ,u_{(0,d_2)} )$, $\bm{u}_i =(u_{(i,0)} ,u_{(i,1)} ,\cdots ,u_{(i,d_2)}) \ (i=1,2, \cdots ,d_1 -1)$, $\widetilde{\bm{u}_{d_1}}=(u_{(d_1 ,0)} ,u_{(d_1 ,1)} ,\cdots ,u_{(d_1 ,d_2 -1)})$. \\
Let us write $(x_1,x_2,\dots, x_{(d_1 +3)(d_2 +3)-6})=(\bm{a}_0 ,\dots ,\bm{a}_{d_1} ,\bm{b}_0 ,\dots ,\bm{b}_{d_2} , \widetilde{\bm{u}_0} ,\bm{u}_1 ,\dots ,\bm{u}_{d_1 -1} ,\widetilde{\bm{u}_{d_1}} )$.  
Then torus actions are given by, 
\begin{align}
&(x_1 ,\dots , x_{(d_1 +3)(d_2 +3)-6}) \notag \\ 
&\ \to \left( \prod _{j=1} ^{(d_1 +1)(d_2 +1)+1} \lambda _j ^{w_{(j,1)}} x_1 ,\dots , \prod _{j=1} ^{(d_1 +1)(d_2 +1)+1} \lambda _j ^{w_{(j,(d_1 +3)(d_2 +3)-6)}} x_{(d_1 +3)(d_2 +3)-6} \right),
\end{align}
where the $((d_1 +1)(d_2 +1)+1) \times ((d_1 +3)(d_2 +3)-6) $ matrix $W_{(d_1 ,d_2)}=(w_{(i,j)})$ is a weight matrix which will be given in Subsection 2.1.
} \label{1.1.2}
\end{defi}

\textit{Ouitline of this paper.} 
\ In Subsection 1.2, we review basic definitions and well-known facts of toric geometry. In Section 2, we explicitly construct toric data of the moduli space $\widetilde{Mp}_{0,2} (\mathbb{P}^1 \times \mathbb{P}^1,\bf{d})$. In Subsection 2.1, we first construct 1--dimensional cones (precisely, generators of $1$--dimensional cones) of $\widetilde{Mp}_{0,2} (\mathbb{P}^1 \times \mathbb{P}^1,\bf{d})$ from its weight matrix 
given in \cite{J1} and fundamental exact sequence used in toric geometry. In Subsection 2.2, we introduce the idea of pseudo-fan and 
construct the fan of $\widetilde{Mp}_{0,2} (\mathbb{P}^1 \times \mathbb{P}^1,\bf{d})$. Then, we reduce the problem to prove that the fan is complete and simplicial, to solving a certain system of linear min--value equations. In Subsection 2.3, we show that the 
definition of $\widetilde{Mp}_{0,2} (\mathbb{P}^1 \times \mathbb{P}^1,\bf{d})$ given in Definition\ 1 follows from the standard quotient construction obtained from the fan given in the previous subsections. In section 2.4, we actually solve the system of min--value equations introduced in Subsection 2.2 and completes the proof of the theorem that $\widetilde{Mp}_{0,2} (\mathbb{P}^1 \times \mathbb{P}^1,\bf{d})$ is a 
compact orbifold. In Section 3, we compute Chow ring of $\widetilde{Mp}_{0,2} (\mathbb{P}^1 \times \mathbb{P}^1,\bf{d})$ and its Poincar\'e polynomial in some lower degree cases.

\subsection{Basic Definitions and Facts in Toric Geometry}

In this subsection, we review standard definitions and basic facts in toric geometry which we will use in this paper. 
For more details, see \cite{CLS} or \cite{F}.                                                                                           
Let $M$ be a free Abelian group of rank $n$ (i.e., $M \simeq \mathbb{Z}^n$) and $N:=\mathrm{Hom} _{\mathbb{Z}} (M ,\mathbb{Z}) (\simeq \mathbb{Z}^n)$ be its dual. We denote $M \otimes \mathbb{R}$ and $N \otimes \mathbb{R}$ by $M_{\mathbb{R}}$ and $N_{\mathbb{R}}$, respsctively. Similarly, we denote that $M_{\mathbb{Q}} := M \otimes \mathbb{Q}$ and $N_{\mathbb{Q}} := N \otimes \mathbb{Q}$. Then we have a natural pairing:
\begin{equation}
\langle \bullet ,\bullet \rangle : M \times N \longrightarrow \mathbb{Z}. \label{1.2.1}
\end{equation}

\begin{defi}
\textit{
\ Let $\sigma$ be a subset of $N_{\mathbb{R}}$.
\begin{enumerate}[(i)]
\item $\sigma$ is a \emph{rational polyhedral cone} if there exists $v_1 , \dots ,v_l \in N$ such that
\begin{equation}
\sigma = \langle v_1 , \dots ,v_l \rangle _{\ge 0} :=\{ \lambda _1 v_1 + \dots + \lambda _l v_l \ \big| \ \lambda _1 , \dots ,\lambda _l \ge 0 \}. \label{1.2.2}
\end{equation}
\item $\sigma$ is a \emph{strongly convex cone} if it is a rational polyhedral cone and satisfies 
\begin{equation}
\sigma \cap (-\sigma ) =\{ 0 \}. \label{1.2.3}
\end{equation}
\item \emph{Dimension} of a cone $\sigma$ is defined by
\begin{equation}
\mathrm{dim}\; \sigma := \mathrm{dim} _{\mathbb{R}}\; \mathrm{span}(\sigma), \label{1.2.4}
\end{equation}
where $\mathrm{span}(\sigma)$ is minimal $\mathbb{R}$--linear subspace including $\sigma$. 
\item \emph{Dual cone} of a cone $\sigma$ is defined by
\begin{equation}
\check{\sigma} := \{ m \in M_{\mathbb{R}} \ \big| \ \langle m ,v \rangle \ge 0 \ \text{for any}\ v \in \sigma \}. \label{1.2.5}
\end{equation}
\item A subset $\tau$ of a cone $\sigma$ is called a \emph{face} of $\sigma$ if there exists $m \in M \cap \check{\sigma}$ such that 
\begin{equation}
\tau = \{ v \in \sigma \ \big| \ \langle m ,v \rangle =0 \}. \label{1.2.6}
\end{equation}
\end{enumerate}
}
\end{defi}

Note that a cone $\sigma :=\langle v_1 ,\dots ,v_l \rangle _{\ge 0}$ generated by $\mathbb{R}$--linearly independent vectors $v_1 ,\dots ,v_l$ is strongly convex. 

Then we introduce basic characteristics of cones without proof. 
\begin{lem}
\textit{
\begin{enumerate}[(i)]
\item A face of a cone is also a cone.
\item Let $\sigma$ be a cone and $\tau$ be its face. Then every face of $\tau$ is a face of $\sigma$.
\item If $\sigma$ is a cone, then $\check{\sigma} \cap M$ is finitely generated semigroup (Gordan's lemma).
\item If $\rho$ is a cone with $\mathrm{dim} \; \rho =1$, then there exists the unique generator $v_{\rho}$ of the semigroup $\rho \cap N$. We identify the $1$--dimensional cone $\rho$ with its generator $v_{\rho}$.
\end{enumerate}
}
\end{lem}
\begin{lem}
\textit{
Let $\{ v_1 , \dots , v_l \}$ be a subset of $N$ and $S$ be a subset of $\{ v_1 , \dots , v_l \}$. If $v_1 , \dots , v_l$ are linearly independent as vectors of $N_{\mathbb{R}}$, then $\langle S \rangle _{\ge 0}$ is a face of $\sigma := \langle v_1 , \dots , v_l \rangle _{\ge 0}$. 
}
\end{lem}

\begin{prf}
\ Since $\sigma$ is a face of $\sigma$ and a face of a face is a face, we can assume that $S=\{ v_1 ,\dots , v_{l-1} \}$. It is clear that $v_1 , \dots , v_l$ are linearly independent as vectors in $N_{\mathbb{Q}}$. We define
\[ V_1 := \{ m \in M_{\mathbb{Q}} \ \big| \ \langle m, v_1 \rangle = \dots = \langle m , v_{l-1} \rangle =0 \}, \]
\[ V_2 := \{ m \in M_{\mathbb{Q}} \ \big| \ \langle m, v_1 \rangle = \dots = \langle m , v_{l} \rangle =0 \}. \]
Then $\mathrm{dim} _{\mathbb{Q}} V_1 = n-(l-1)$ and $\mathrm{dim} _{\mathbb{Q}} V_2 = n-l$. Thus $V_1 \supsetneq V_2$ and we can take $m_0 \in V_1 \setminus V_2$. By  multiplying $m_0$ by some integer if necessary, $m_0$ can be made into an element of  $M$
that satisfies $\langle m_0 , v_l \rangle >0$. This means that 
\[ \langle S \rangle _{\ge 0} = \{ v \in \sigma \ \big| \ \langle m_0 ,v \rangle =0 \}, \]
and therefore $\langle S \rangle _{\ge 0}$ is a face of $\sigma$. \ $\square$
\end{prf}

\begin{defi}
\textit{\ A set $\Sigma$ which consists of strongly convex cones is called a \emph{fan} if
\begin{enumerate}[(i)]
\item $\Sigma$ is a finite set.
\item for any $\sigma \in \Sigma$, every face of $\sigma$ is also cone of $\Sigma$.
\item for any $\sigma , \tau \in \Sigma$, $\sigma \cap \tau$ is a face of $\sigma$ and $\tau$.
\end{enumerate}
Moreover, the set $|\Sigma|$ denotes $\bigcup _{\sigma \in \Sigma} \sigma$ and we define the set
\begin{equation}
\Sigma (k) := \{ \sigma \in \Sigma \ \big| \ \mathrm{dim}\; \sigma =k \}. \label{1.2.7}
\end{equation}
Especially, we also denote by $\Sigma(1)$ a set of generators of 1-dimensional cones in $\Sigma$. 
}
\end{defi}

A fan $\Sigma$ defines a toric variety $X(\Sigma)$ which is obtained by gluing \textit{affine toric varieties} $X_{\sigma} :=\mathrm{Spec}(\mathbb{C} [M \cap \check{\sigma}])\ (\sigma \in \Sigma)$. Then $X_{\{ 0 \}} =\mathrm{Spec}(\mathbb{C} [M])$ is an algebraic torus $T_N :=N \otimes \mathbb{C}^{\times} \simeq (\mathbb{C}^{\times})^n$.

At this stage, we introduce a fundamental theorem in toric geometry which we use in this paper. First, we give the geometrical properties of $X(\Sigma)$ using cones of $\Sigma$:   
\begin{thm}
\textit{\ Let $\Sigma$ be a fan.
\begin{enumerate}[(i)]  
\item $X(\Sigma)$ is simplicial (i.e., an orbifold) if and only if $\Sigma$ is \emph{simplicial} (i.e., the generators of $\sigma$ are linearly independent as vectors of $N_{\mathbb{R}}$ for any $\sigma \in \Sigma$).
\item $X(\Sigma)$ is complete (i.e., compact) if and only if $\Sigma$ is \emph{complete} (i.e., $|\Sigma|=N_{\mathbb{R}}$).
\end{enumerate}
}
\end{thm}

If the generators of $1$--dimentional cones of $\Sigma$ span $N_{\mathbb{R}}$, then the following sequence of $\mathbb{Z}$--modules is exact:
\begin{equation}
\xymatrix@C=40pt{
0 \ar[r] & M \ar[r]^-{\varphi} & \mathbb{Z} ^{\Sigma (1)} \ar[r]^-{\psi} & A_{n-1} (X(\Sigma)) \ar[r] & 0 
}, \label{1.2.8}
\end{equation}
where $\varphi (m) := ( \langle m ,v_{\rho}  \rangle )_{\rho}$, and $\psi ((a_{\rho})_{\rho})$ is the divisor class of $\sum _{\rho} a_{\rho} D_{\rho} \, $ (in general, $1$--dimensional cone $\rho$ of $\Sigma$ corresponds to $T_{N}$--invariant divisor $D_{\rho}$ on $X(\Sigma)$).

As is well known in the case of complex projective space, we can obtain the quotient construction (i.e., representation using homogeneous coordinates) of toric variety $X(\Sigma)$ when the fan $\Sigma$ is complete and simplicial (see \cite{C1}). 
We briefly explain outline of this construction. 
First, we define a closed subspace of $\mathbb{C}^{\Sigma (1)}$ by,
\begin{equation}
Z(\Sigma) :=\left\{ (x_{\rho})_{\rho \in \Sigma (1)} \in \mathbb{C}^{\Sigma (1)} \ \Bigg| \ \prod_{\rho \not\subset \sigma} x_{\rho} =0 \ \text{for all} \ \sigma \in \Sigma \right\} \label{1.2.9}
\end{equation}
and an abelian group that acts on $\mathbb{C}^{\Sigma (1)}$ by,
\begin{equation}
G:=\mathrm{Hom} _{\mathbb{Z}} (A_{n-1} (X(\Sigma)),\mathbb{C}^{\times}), \label{1.2.10}
\end{equation}
where $A_{n-1} (X(\Sigma))$ is Chow group of $X(\Sigma)$ given in (\ref{1.2.8}). The group action is given by,
\begin{equation}
G \times \mathbb{C}^{\Sigma (1)} \ni (g, (x_{\rho})_{\rho}) \longmapsto (g([D_{\rho}])x_{\rho})_{\rho} \in \mathbb{C}^{\Sigma (1)}. \label{1.2.11}
\end{equation}
Then we can state the following theorem from \cite{C1}:
\begin{thm}
\textit{\ Let $\Sigma$ be a fan that the generators of $1$--dimentional cones of $\Sigma$ span $N_{\mathbb{R}}$. Then $X(\Sigma)$ is geometrical quotient of $\mathbb{C}^{\Sigma (1)} - Z(\Sigma)$ by $G$ if and only if $\Sigma$ is simplicial.}
\end{thm}

In order to describe the closed set $Z(\Sigma)$ in a simpler manner, we introduce \textit{primitive collection}. In the following, we 
assume that $\Sigma$ is simplicial. 
\begin{defi}(\cite{BC})                                                                                                                                                                   
\textit{\ Let $\Sigma$ be a simplicial fan. A subset $S$ of $\Sigma(1)$ is a} primitive collection \textit{of $\Sigma$ if 
\begin{enumerate}
\item[(i)]
$S$ does not generate a cone in $\Sigma$.
\item[(ii)] 
every proper subset of $S$ generate a cone in $\Sigma$.
\end{enumerate}
Moreover, we denote set of primitive collections of $\Sigma$ by $\mathrm{PC} (\Sigma)$.
}
\end{defi} 

If $\Sigma$ is simplicial, $Z(\Sigma)$ can be expressed in the following form by using the primitive collections:
\begin{lem}(\cite{BC})                                                                                                                                                          
\textit{\ If $\Sigma$ is simplicial, then
\begin{equation}
Z(\Sigma) =\bigcup _{S \in \mathrm{PC} (\Sigma)}  \{ (x_{\rho}) \in \mathbb{C} ^{\Sigma (1)} \ \big| \ x_{\rho} =0 \ \text{for any} \ \rho \in S \}.
\end{equation}
}
\end{lem}

\textit{Chow Ring and Cohomology Ring of Toric Varieties.}
\ \ Let $\Sigma$ be an $n$--dimensional fan. In order to compute Chow Ring of $X(\Sigma)$, we define two ideals of polynomial ring $\mathbb{Q}[x_{\rho} | \rho \in \Sigma (1)]$. The first ideal which we denote by $I(\Sigma)$, is defined as an ideal generated by homogeneous polynomials $\sum_{\rho \in \Sigma(1)} \langle m,v_{\rho} \rangle \cdot x_{\rho} \ (m \in M)$. 
By taking $m$ as canonical basis $e_{i}\;\;(i=1,\cdots,n)$ of $M$, this ideal is equal to
\begin{equation}
\left\{ \sum _{\rho \in \Sigma (1)} v_{i \rho} x_{\rho} \Bigg| i=1 , \dots , n \right\},
\end{equation}
where $v_{\rho} ={}^t\hspace{0.5pt}(v_{1 \rho} , \dots , v_{n \rho})$. The second ideal, called \textit{Stanley--Reisner ideal}, is defined by,
\begin{equation}
SR(\Sigma) := \left\{ \prod _{v_{\rho} \in S} x_{\rho} \Bigg| S \in \mathrm{PC}(\Sigma) \right\}.
\end{equation}
Then the Chow ring, (roughly speaking, cohomology ring) , and Betti numbers of complete simplicial toric variety $X(\Sigma)$ can be computed as follows \cite{CK,CLS}:
\begin{thm}
\textit{\ Let $\Sigma$ be complete simplicial fan. Then
\begin{enumerate}[(i)]
\item 
\begin{equation}
A_{\mathbb{Q}} ^{\ast} (X(\Sigma)) \simeq H^{\ast} (X(\Sigma) ,\mathbb{Q}) \simeq \mathbb{Q}[x_{\rho} | \rho \in \Sigma (1)]/(I(\Sigma) +SR(\Sigma)).  \label{betti}
\end{equation}
\item
\begin{equation}
b_{2k} (X(\Sigma)) = \sum _{i=k}^{n} (-1)^{i-k} \dbinom{i}{k} |\Sigma(n-i)| \quad (\text{for all $k$}),
\end{equation}
where $b_{2k} (X(\Sigma)) :=\mathrm{dim} \, H^{2k} (X(\Sigma),\mathbb{Q})$ is Betti number of $X(\Sigma)$.
\item
\begin{equation}
b_{2k+1} (X(\Sigma))=0 \quad (\text{for all $k$}).
\end{equation}
\item
\begin{equation}
b_{2k} (X(\Sigma)) = b_{2n-2k} (X(\Sigma)) \quad (\text{for all $k$}).
\end{equation}
\end{enumerate}
}
\end{thm}
Using this theorem, we can also compute Poincar\'e Polynomial $P(t):=\sum_{i=0}^{2n} b_{2i} t^i$ of complete simplicial toric variety by 
counting the number of cones in its fan.

\vspace{2em}
\noindent
{\bf Acknowledgment.} \ The author would like to thank Prof. M. Jinzenji for suggesting the problem treated in this paper, his big support and many useful discussions.

\section{Toric Construction of $\widetilde{Mp}_{0,2} (\mathbb{P}^1 \times \mathbb{P}^1 , (d_1 ,d_2))$}

\subsection{Weight and Vertex Matrix}

In order to define weight and vertex matrices of $\widetilde{Mp}_{0,2} (\mathbb{P}^1 \times \mathbb{P}^1 , (d_1 ,d_2))$, we label rows and columns of these matrices by using the following column and row vectors:

{\bf column vectors used for labeling row vectors of matrices} :
\[ z={}^t\hspace{0.5pt}(z_0 ,z_1 , \dots ,z_{d_1}), \]
\[ w={}^t\hspace{0.5pt}(w_0 ,w_1 , \dots ,w_{d_2}), \]
\[ f_i ={}^t\hspace{0.5pt}(f_{(i,1)} ,f_{(i,2)} ,\dots ,f_{(i,d_2)} ) \ (i=1 ,2 ,\dots ,d_1), \]
\[ G_{d_1} ={}^t\hspace{0.5pt}(g_1 ^{(d_1 )} , g_2 ^{(d_1 )} , \dots ,g_{d_1 -1} ^{(d_1 )}), \]
\[ G_{d_2} ={}^t\hspace{0.5pt}(g_1 ^{(d_2 )} , g_2 ^{(d_2 )} , \dots ,g_{d_2 -1} ^{(d_2 )}), \]
\[ G ={}^t\hspace{0.5pt}(g), \]

{\bf row vectors used for labeling column vectors of matrices} :
\[ a^i =(a_0 ^i ,a_1 ^i ,\dots ,a_{d_1} ^i) \ (i=1,2), \]
\[ b^i =(b_0 ^i ,b_1 ^i ,\dots ,b_{d_2} ^i) \ (i=1,2), \]
\[ \widetilde{u_0} =(u_{(0,1)} ,u_{(0,2)}, \cdots ,u_{(0,d_2)} ), \]
\[ u_i =(u_{(i,0)} ,u_{(i,1)} ,\dots ,u_{(i,d_2)}) \ (i=1,2, \dots ,d_1 -1), \]
\[ \widetilde{u_{d_1}}=(u_{(d_1 ,0)} ,u_{(d_1 ,1)} ,\dots ,u_{(d_1 ,d_2 -1)}). \]
We introduce the following two sets which are used in subsection 2.2: 
\[ \mathrm{Row} (d_1 ,d_2 ):=\{ z_0 ,\dots ,z_{d_1} ,w_0 ,\dots ,w_{d_2} ,f_{(1,1)}, \dots ,f_{(d_1 ,d_2)} \}, \] 
and 
\[ \mathrm{Col} (d_1 ,d_2) :=\{ a_0 , \dots a_{d_1} , b_0 ,\dots ,b_{d_2} ,u_{(0,1)} ,\dots u_{(0,d_2)}, u_{(1,0)}, \dots , u_{(1,d_2)}, \dots ,u_{(d_1 ,0)}, \dots ,u_{(d_1 ,d_2 -1)} \}, \] 
where $a_i := (a_i ^1 ,a_i ^2) , b_j := (b_j ^1 ,b_j ^2) \ (i=0, \dots ,d_1 ,j=0 ,\dots , d_2)$. 
We denote the canonical $\mathbb{Z}$--bases $e_i$, $e_{d_1 +1 +j}$, and $e_{d_1 + d_2 +2 +(i-1)d_2 +j}$ of $\mathbb{Z}^{(d_1 +1)(d_2 +1) +1}$ by $e_{z_i}$, $e_{w_j}$, and $e_{f_{(i,j)}}$, respsctively. Moreover, we put $n:=2d_1 +2d_2 +1$ and $r:=|\Sigma_{(d_1 ,d_2 )} (1)|=(d_1 +3)(d_2 +3)-6$. Then $r-n=(d_1 +1)(d_2 +1)+1$. With this set up, we define the $(r-n)\times r$ weight matrix $W_{(d_1 ,d_2)}$\ ($d_1 \ge d_2 >0$) by,

\begin{align}
&W_{(d_1 ,d_2)} = (w_{(i,j)}) \notag \\
&:= \bordermatrix{
                         & a^1 & a^2 & b^1 & b^2 & \widetilde{u_0} & u_1 & u_2 & \cdots & u_{d_1 -2} & u_{d_1 -1} & \widetilde{u_{d_1}} \cr 
                       z & I_{d_1 +1} & I_{d_1 +1} & O & O & U_0 & K_1 & K_2 & \cdots &K_{d_1 -2} & K_{d_1 -1} & U_{d_1} \cr
                       w & O & O & I_{d_2 +1} & I_{d_2 +1} & L_0 & O & O & \cdots & O & O & L_{d_1} \cr
                       f_1 & O & O & O & O & \widetilde{J_d} & J_u & O & \cdots & O & O & O \cr
                       f_2 & O & O & O & O & O & J_d & J_u & \cdots & O & O & O \cr
                       \vdots & \vdots & \vdots & \vdots & \vdots & \vdots & \vdots & \vdots & \ddots & \vdots & \vdots & \vdots \cr
                       f_{d_1 -1} & O & O & O & O & O & O & O & \cdots & J_d & J_u & O \cr
                       f_{d_1} & O & O & O & O & O & O & O & \cdots & O & J_d & \widetilde{J_u}  
                      },
\end{align}
where $O$ is zero matrix and $I_N \ (N \in \mathbb{N})$ is identity matrix of size $N$, and other matrices are defined as follows:
\begin{equation}
U_0 := \left( 
\begin{array}{ccc|c}
         &&& 0 \\
         &&& -1 \\
         & \text{\huge{0}} &&0 \\
         &&& \vdots \\
         &&& 0 
         \end{array}
         \right) ,
\end{equation}
\begin{equation}
U_{d_1} := \left(
\begin{array}{c|ccc}
        0 &&& \\
        \vdots &&& \\
        0 && \text{\huge{0}} & \\
        -1 &&& \\
        0 &&& 
        \end{array}
        \right) ,
\end{equation}
\begin{equation}
K_i :=
\begin{array}{rc|ccc|cll}
        \ldelim({10}{2mm} & 0 &&&& 0 & \rdelim){10}{2mm} & \\ 
       & \vdots &&&& \vdots && \\
       & 0 &&&& 0 && \\
       & -1 &&&& 0 && \\
       & 1 && \text{\huge{0}} && 1 && \!\!\!<z_i \\
       & 0 &&&& -1 && \\
       & 0 &&&& 0 && \\
       & \vdots &&&& \vdots && \\
       & 0 &&&& 0 &&
        \end{array} 
         \quad (i=1,\dots , d_1 -1) ,
\end{equation}
\begin{equation}
L_0 := \left(
\begin{array}{cccccc|c}
-1 &&&&&& 0 \\
1 & -1 &&& \text{\huge{0}} && 0 \\
& 1 & -1 &&&& 0 \\
&&& \ddots &&& \vdots \\
& \text{\huge{0}} &&& 1 & -1 & 0 \\
&&&&& 1 & -1 \\ \hline
0 & 0 & 0 & \ldots & 0 & 0 & 0
\end{array}
\right) ,
\end{equation}
\begin{equation}
L_{d_1} := \left(
\begin{array}{c|cccccc}
0 & 0 & 0 & \ldots & 0 & 0 & 0 \\ \hline
-1 & 1 &&&&& \\
0 & -1 & 1 &&& \text{\huge{0}} & \\
\vdots &&& \ddots &&& \\
0 &&&  & -1 & 1 & \\
0 && \text{\huge{0}} &&& -1 & 1 \\
0 &&&&&& -1
\end{array} 
\right) ,
\end{equation}
\begin{equation}
\widetilde{J_d} := \left(
\begin{array}{cccccc}
1 &&&&& \\
-1 & 1 &&& \text{\huge{0}} & \\
& -1 & 1 &&& \\
&&& \ddots && \\
& \text{\huge{0}} &&& 1 & \\
&&&& -1 & 1
\end{array}
\right) ,
\end{equation}
\begin{equation}
\widetilde{J_u} := \left(
\begin{array}{cccccc}
1 & -1 &&&& \\
& 1 & -1 && \text{\huge{0}} & \\
&& 1 &&& \\
&&& \ddots && \\
& \text{\huge{0}} &&& 1 & -1 \\
&&&&& 1
\end{array}
\right) ,
\end{equation}
\begin{equation}
J_d := \left(
\begin{array}{cccccc}
-1 &1 &&&& \\
& -1 & 1 && \text{\huge{0}} & \\
&& -1 &&& \\
&&& \ddots && \\
& \text{\huge{0}} &&& 1 & \\
&&&& -1 & 1
\end{array}
\right) ,
\end{equation}
\begin{equation}
\widetilde{J_u} := \left(
\begin{array}{cccccc}
1 & -1 &&&& \\
& 1 & -1 && \text{\huge{0}} & \\
&& 1 &&& \\
&&& \ddots && \\
& \text{\huge{0}} &&& -1 & \\
&&&& 1 & -1
\end{array}
\right).
\end{equation}
In some lower degree cases, this weight matrix was presented in \cite{J1} by using degree diagram, but we give here explicit definition 
of general $W_{(d_1,d_2)}$ by generalizing the construction in \cite{J1}.

In order to define the fan of $\widetilde{Mp}_{0,2} (\mathbb{P}^1 \times \mathbb{P}^1 , (d_1 ,d_2))$, we first introduce the matrix that gives generators of $1$--dimensional cones of the fan.  In this paper, we call this matrix \textit{vertex matrix}. For example, the fan of $2$--dimensional projective space $\mathbb{P}^2$ consists of 
$\{0\},
 \langle v_1 \rangle _{\ge 0},
 \langle v_2 \rangle _{\ge 0},
 \langle v_3 \rangle _{\ge 0},$
$ 
 \langle v_1 ,v_2 \rangle _{\ge 0},
 \langle v_1 ,v_3 \rangle _{\ge 0},$
  and
$  
 \langle v_2 ,v_3 \rangle _{\ge 0},$
 where $v_1 := {}^t\hspace{0.5pt}(1,0), v_2 :={}^t\hspace{0.5pt}(0,1), v_3 :={}^t\hspace{0.5pt}(-1,-1)$.
Hence the vertex matrix of $\mathbb{P}^2$ is given by,
\begin{equation}
(v_1 ,v_2 ,v_3 ) = \left( \begin{array}{ccc}
        1 & 0 & -1 \\
        0 & 1 & -1
        \end{array} \right).
\end{equation}
Then we define the $n\times r$ vertex matrix $V_{(d_1 ,d_2 )}$ of $\widetilde{Mp}_{0,2} (\mathbb{P}^1 \times \mathbb{P}^1 , (d_1 ,d_2))$ as follows:
\begin{align}
&V_{(d_1 ,d_2)} = (v_1 , \dots ,v_r ) \ \notag \\ 
&:= \bordermatrix{
                         & a^1 & a^2 & b^1 & b^2 & \widetilde{u_0} & u_1 & u_2 & \cdots & u_{d_1 -2} & u_{d_1 -1} & \widetilde{u_{d_1}} \cr 
                       z & I_{d_1 +1} & -I_{d_1 +1} & O & O & O & O & O & \cdots & O & O & O \cr
                       w & O & O & I_{d_2 +1} & -I_{d_2 +1} & O & O & O & \cdots & O & O & O \cr
                       G_{d_1} & O & \widetilde{A_{d_1}} & O & O & O & \widetilde{e_1} & \widetilde{e_2} & \cdots & \widetilde{e_{d_1 -2}} & \widetilde{e_{d_1 -1}} & O \cr
                       G & O & G_{a} & O & G_{b} & -1^{d_2} & e_1 & e_1 & \cdots & e_1 & e_1 & e_1 ^{-} \cr
                       G_{d_2} & O & O & O & \widetilde{A_{d_2}} & I_{d_2} ^{R} & I_{d_2} ^{LR} & I_{d_2} ^{LR} & \cdots & I_{d_2} ^{LR} & I_{d_2} ^{LR} & I_{d_2} ^{L} \cr
                      }, \label{vertex}
\end{align}
where $e_1$ and $e_1 ^{-}$ are first canonical bases of $\mathbb{Z} ^{d_2 +1}$ and $\mathbb{Z} ^{d_2}$, respectivery, and
\begin{equation}
\widetilde{A_{N}} :=\left(
\begin{array}{ccccccccc}
1 & -2 & 1 &&&&&& \\
& 1 & -2 & 1 &&&& \text{\huge{0}}& \\
&& 1 & -2 & 1 &&&& \\ 
& \text{\huge{0}} &&&& \ddots &&& \\
&&&&&& 1 & -2 & 1
\end{array}
\right)  \in M(N-1,N+1;\mathbb{Z}),
\end{equation}
\begin{equation}
G_a := (1 ,-1,0, \dots ,0),
\end{equation}
\begin{equation}
G_b := (-1 ,1,0, \dots ,0),
\end{equation}
\begin{equation}
-1^{d_2} := (-1 ,-1, \dots ,-1),
\end{equation}
\begin{equation}
\widetilde{e_i} := 
\begin{array}{rccccll}
\ldelim({8}{2mm}&0&0&\ldots&0 &\rdelim){8}{2mm} & \\
&\vdots&\vdots&&\vdots && \\
&0&0&\ldots&0 && \\
&1&1&\ldots&1 && \!\!\!<g_i ^{(d_1)} \\
&0&0&\ldots&0 && \\
&\vdots&\vdots&&\vdots && \\
&0&0&\ldots&0 && 
\end{array}
\ (i=1, \dots ,d_1),
\end{equation}
\begin{equation}
I_{d_2} ^{R} := \left(
\begin{array}{ccccc|c}
1&&&&&0 \\
&1&&\text{\huge{0}}&&0 \\
&&\ddots&&&\vdots \\
&\text{\huge{0}}&&1&&0 \\
&&&&1&0
\end{array}
\right) ,
\end{equation}
\begin{equation}
I_{d_2} ^{LR} := \left(
\begin{array}{c|ccccc|c}
0&1&&&&&0 \\
0&&1&&\text{\huge{0}}&&0 \\
\vdots&&&\ddots&&&\vdots \\
0&&\text{\huge{0}}&&1&&0 \\
0&&&&&1&0
\end{array}
\right) ,
\end{equation}
\begin{equation}
I_{d_2} ^{L} := \left(
\begin{array}{c|ccccc}
0&1&&&& \\
0&&1&&\text{\huge{0}}& \\
\vdots&&&\ddots&& \\
0&&\text{\huge{0}}&&1& \\
0&&&&&1
\end{array}
\right).
\end{equation}

\begin{ex}
\begin{equation}
W_{(1,1)} =\left( \begin{array}{cccccccccc}
1 & 0 & 1 & 0 & 0 & 0 & 0 & 0 & 0 & -1 \\
0 & 1 & 0 & 1 & 0 & 0 & 0 & 0 & -1 & 0 \\ 
0 & 0 & 0 & 0 & 1 & 0 & 1 & 0 & -1 & 0 \\
0 & 0 & 0 & 0 & 0 & 1 & 0 & 1 & 0 & -1 \\
0 & 0 & 0 & 0 & 0 & 0 & 0 & 0 & 1 & 1 
\end{array} \right)
\end{equation}
\begin{equation}
V_{(1,1)} =\left( \begin{array}{cccccccccc}
1 & 0 & -1 & 0 & 0 & 0 & 0 & 0 & 0 & 0 \\
0 & 1 & 0 & -1 & 0 & 0 & 0 & 0 & 0 & 0 \\
0 & 0 & 0 & 0 & 1 & 0 & -1 & 0 & 0 & 0 \\
0 & 0 & 0 & 0 & 0 & 1 & 0 & -1 & 0 & 0 \\
0 & 0 & 1 & -1 & 0 & 0 & -1 & 1 & -1 & 1
\end{array} \right)
\end{equation}
\end{ex}

The following lemmas are easily confirmed by using the definitions of $W_{(d_1 ,d_2)}$ and $V_{(d_1 ,d_2)}$.

\begin{lem}
\textit{\ The sequence
\begin{equation}
\xymatrix@C=48pt{
0 \ar[r] & \mathbb{Z}^{n} \ar[r]^-{{}^t\hspace{0.5pt} V_{(d_1 ,d_2)}} & \mathbb{Z}^{r} \ar[r]^-{W_{(d_1 ,d_2)}} & \mathbb{Z}^{r-n} \ar[r] & 0
} \label{foundexactmat}
\end{equation}
is exact.
}
\end{lem}

\begin{lem}
\textit{
\begin{enumerate}[(i)]
\item
\begin{equation} 
\sum _{j=0} ^{d_2} \mathbf{w}_{u_{(i,j)}} = -e_{z_{i -1}} +2e_{z_i} -e_{z_{i+1}} \ (i=1 ,\dots ,d_1 -1), \label{2.1.24}
\end{equation}
\item
\begin{equation} 
\sum _{i=0} ^{d_1} \mathbf{w}_{u_{(i,j)}} = -e_{w_{j -1}} +2e_{w_j} -e_{w_{j+1}} \ (j=1 ,\dots ,d_2 -1), \label{2.1.25}
\end{equation}
\item
\begin{equation}
\left( \sum _{p=0} ^{i-1} \sum _{q=j+1} ^{d_2} \mathbf{w} _{u_{(p,q)}} \right)_{\!\!f_{(s,t)}} = \left( e_{f_{(i ,j+1)}} \right)_{f_{(s,t)}} \ (s=1, \dots ,d_1 ,t=1 \dots ,d_2), \label{2.1.26}
\end{equation}
\item
\begin{equation}
\left( \sum _{p=i+1} ^{d_1} \sum _{q=0} ^{j-1} \mathbf{w} _{u_{(p,q)}} \right)_{\!\!f_{(s,t)}} = \left( e_{f_{(i+1 ,j)}} \right)_{f_{(s,t)}} \ (s=1, \dots ,d_1 ,t=1 \dots ,d_2), \label{2.1.27}
\end{equation}
\end{enumerate}
where $\mathbf{w}_{u_{(i,j)}}$ is the $u_{(i,j)}$--columun vector of $W_{(d_1 ,d_2)}$ and $( \mathbf{x} )_{f_{(s,t)}}$ is the $f_{(s,t)}$--component of a vector
\begin{equation}
\mathbf{x}=(x_{z_0} ,\dots , x_{z_{d_1}} ,x_{w_0} \dots ,x_{w_{d_2}} ,x_{f_{(1,1,)}} ,\dots ,x_{f_{(1,d_2 )}},\dots ,x_{f_{(d_1 ,1)}} ,\dots ,x_{f{(d_1 ,d_2)}}). \label{2.1.28}
\end{equation}
}
\end{lem}

\subsection{Pseudo-Fan and Construction of Fan}

We denote the $\rho$-th column vector of $V_{(d_1 ,d_2)}$ by $v_{\rho}$ and a collection of column vectors of $V_{(d_1 ,d_2)}$ by $\mathrm{Ver} _{(d_1 ,d_2)}$. In this subsection, we construct a fan that leads us to the definition of $\widetilde{Mp}_{0,2} (\mathbb{P}^1 \times \mathbb{P}^1 , (d_1 ,d_2))$ given in Definition 1.
Let us illustrate our strategy of construction by taking the $(d_1 ,d_2 )=(1,1)$ case as an example. In this case,
 $\widetilde{Mp}_{0,2} (\mathbb{P}^1 \times \mathbb{P}^1 , (1,1))$ can be written as,
\begin{align}
& \{ ( \bm{a}_0 ,\bm{a}_1 ,\bm{b}_0 ,\bm{b}_1 ,u_{(0,1)} ,u_{(1,0)}) \in \mathbb{C}^{10} \big| \bm{a}_0 ,\bm{a}_{d_1},\bm{b}_0 ,\bm{b}_{d_2} \neq \bm{0}, (u_{(0,1)} 
,u_{(1,0)})\neq (0,0) \} / (\mathbb{C}^{\times})^{5} \notag \\
& = (\mathbb{C}^{10} - (\{ \bm{a}_0 = \bm{0} \} \cup \{ \bm{a}_1 =\bm{0} \} \cup \{ \bm{b}_0 = \bm{0} \} \cup \{ \bm{b}_1 =\bm{0} \} \cup \{ (u_{(0,1)} 
,u_{(1,0)}) = (0,0) \}))/(\mathbb{C}^{\times})^{5}.
\end{align}
Therefore, by Theorem 2, we should define a fan $\Sigma _{(1,1)}$ of $\widetilde{Mp}_{0,2} (\mathbb{P}^1 \times \mathbb{P}^1 , (1,1))$ so that the following equality holds: 
\begin{equation}
Z(\Sigma_{(1,1)}) = \{ \bm{a}_0 = \bm{0} \} \cup \{ \bm{a}_1 =\bm{0} \} \cup \{ \bm{b}_0 = \bm{0} \} \cup \{ \bm{b}_1 =\bm{0} \} \cup \{ (u_{(0,1)} 
,u_{(1,0)}) = (0,0) \}.
\end{equation}

With Lemma 3 in mind, we construct $\Sigma _{(1,1)}$ as a collection of cones generated by subsets of $\mathrm{Ver} _{(1,1)}$ such that every member of the subsets does not contain $\{v_{a_0^1} ,v_{a_0^2}\}$, $\{v_{a_1 ^1} ,v_{a_1 ^2}\}$, $\{v_{b_0 ^1} ,v_{b_0 ^2}\}$, $\{ v_{b_1 ^1} ,v_{b_1 ^2}\}$, and $\{ v_{u_{(0,1)}},v_{u_{(1,0)}}\}$ (see Lemma 6 below).

In order to define the fan in general case, we interpret the fan from the point of view of primitive collections:
\begin{defi}
\textit{\ Let $V$ be a finite subset of $N_{\mathbb{R}} \simeq \mathbb{R}^d$ and $P_1 , \dots ,P_l$ be non-empty subsets of $V$. Moreover, we put $\Pi := \{ P_1 , \dots , P_l \}$ and
\begin{equation}
\Sigma = \Sigma (V,\Pi ) := \{ \langle S \rangle _{\ge 0} \ \big| \ S \subset V \text{, and $P_1 ,\dots ,P_l$ are not contained in $S$ } \}. \label{2.3.1}
\end{equation}
In this paper, we call a triplet $(\Sigma , V, \Pi)$} $d$--dimensional pseudo-fan \textit{if
\begin{enumerate}[(i)]
\item
\begin{equation}
P_1 \cup \dots \cup P_l = V. \label{2.3.2}
\end{equation}
\item
Any $P_i$ is not contained in another element of $\Pi$.
\end{enumerate}
}
\end{defi}                                                                                           
The following condition (MVC) gives a sufficient condition for a pseudo fan to be a complete and simplical fan.
\begin{defi}
\textit{\ Let $(\Sigma ,V, \Pi)$ be a $d$--dimensional pseudo-fan. Then we say that $(\Sigma ,V, \Pi)$ satisfies \emph{Min--Value Condition (MVC)} if for any $\beta \in \mathbb{R}^d$, there exists the unique $\alpha =(\alpha _i) \in \mathbb{R}^{r} \ (r:=|V|)$ such that
\begin{enumerate}[(i)]
\item 
\begin{equation}
\beta = \sum _{i=1} ^{r} \alpha _i v_i, \label{2.3.3}
\end{equation} 
where $V=\{ v_1 ,\dots , v_r \}$.
\item 
\begin{equation}
\mathrm{min} \{ \alpha _i \ | \ v_i \in P \} =0 \ \ (\text{for any}\ P \in \Pi). \label{2.3.4}
\end{equation}
\end{enumerate}
}
\end{defi}

\begin{lem}
\textit{\ Let $(\Sigma ,V ,\Pi)$ be a $d$--dimensional pseudo-fan satisfying (MVC). Then $\Sigma$ is a $d$--dimentional complete and simplicial fan with $\mathrm{PC}(\Sigma) =\Pi$}.
\end{lem}

\begin{prf}
Let $\sigma := \langle u_1 ,\dots ,u_l \rangle _{\ge 0} \ (u_1 ,\dots ,u_l \in V)$ be any member of $\Sigma$. If $u_1 ,\dots ,u_l$ are linearly dependent over $\mathbb{R}$, there exists $\alpha _1 , \dots ,\alpha _l \in \mathbb{R}$ such that $\alpha _1 u_1 + \dots + \alpha _l u_l =0$, $(\alpha _1 ,\dots ,\alpha _l ) \neq (0, \dots ,0)$, and we can assume that $\alpha _1 ,\dots ,\alpha _p \ge 0 $, $\alpha _{p+1} ,\dots ,\alpha _{l} <0 $ ($1 \le p \le l$). Then there exists $v \in V \setminus \{ u_1 ,\dots ,u_l \}$ such that $v \in P$ for each $P \in \Pi$ since every $P \in \Pi$ is not contained in $\{ u_1 ,\dots ,u_l \}$ (see (\ref{2.3.1})), and therefore both representations 
\begin{equation}
\alpha _1 u_1 + \dots +\alpha _p u_p +\sum_{v \in V \setminus \{ u_1 ,\dots ,u_p \}} 0 \cdot v =(-\alpha _{p+1}) u_{p+1} + \dots +(-\alpha _l )u_l + \sum_{v \in V \setminus \{ u_{p+1} ,\dots ,u_l \}} 0 \cdot v \label{2.3.5}
\end{equation}
satisfy (i),(ii) of (MVC). This follows that $(\alpha _1 ,\dots , \alpha _l )=(0,\dots ,0)$ by uniqueness of $\alpha$, which contradicts our assumption. Thus, $u_1 ,\dots ,u_l$ are linearly independent. Hence $\Sigma$ is simplicial. In particular, every cone in $\Sigma$ is strongly convex.

Let $\beta$ be any member of $N_{\mathbb{R}} \simeq \mathbb{R}^d$. Then we can uniquely write as $\beta = \sum_{i=1}^l \alpha_i v_i$ ($V = \{ v_1 , \dots , v_r \}$) and $\mathrm{min} \{ \alpha_i \ | \ v_i \in P \} = 0$ for any $P \in \Pi$. We set $S^{+} := \{ i \ | \ \alpha_i >0 \}$ and $S^0 := \{ i \ | \ \alpha_i = 0 \}$. It is clear that $S^{+} \cap S^{0} = \emptyset$, and from the condition (i) in Definition 5, $S^{+} \cup S^{0} = V$. Since we can take $v_i \in P$ such that $i \in S^0$ for each $P \in \Pi$, all members of $\Pi$ are not contained in the set $S := \{ v_i \ | \ i \in S^{+} \}$. This means that $\langle S \rangle _{\ge 0} \in \Sigma$, and then $\beta = \sum_{i \in S^{+}} \alpha_i v_i \in \langle S \rangle _{\ge 0}$. Thus $| \Sigma | = N_{\mathbb{R}}$. 

Next, let $\sigma := \langle S \rangle _{\ge 0}$, $\tau := \langle T \rangle _{\ge 0} \ (S,T \subset V)$ be any members of $\Sigma$, and we put $S \cap T=\{ u_1 , \dots ,u_l \}$, $S \setminus T =\{ w_1 ,\dots ,w_p \}$, and $T \setminus S =\{ w_1 ^{\prime} ,\dots ,w_q ^{\prime} \}$. It is clear that $\langle S \cap T \rangle _{\ge 0} \subset \sigma \cap \tau$. Conversely, if $u \in \sigma \cap \tau$, we can write 
\begin{equation}
u=\sum _{i=1} ^{l} \alpha _i u_i + \sum _{j=1} ^{p} \gamma _j w_j +\sum_{v \in S} 0 \cdot v = \sum _{i=1} ^{l} \alpha _i ^{\prime} u_i + \sum _{k=1} ^{q} \gamma _k ^{\prime} w_k ^{\prime} + \sum_{v \in T} 0 \cdot v \label{2.3.6}
\end{equation}
for some $\alpha _i ,\alpha _i ^{\prime} ,\gamma _j , \gamma _k ^{\prime} \ge 0$.
Since these representations satisfy (i),(ii) of (MVC) as above, $\alpha _i = \alpha _i ^{\prime}$ and $\gamma _j , \gamma _k ^{\prime} =0$ by uniqueness. This means that $\sigma \cap \tau = \langle S \cap T \rangle _{\ge 0}$ and $\Sigma$ satisfies (iii) of Definition 3 by Lemma 2. The conditions (i), (ii) of Definition 3 are obvious by (\ref{2.3.1}) and Lemma 2, respectively (as for (ii) of Definition 3, if $P \in \Pi$ is not contained in S, every subset of $S$ does not contain $P$). Therefore, $\Sigma$ is  a complete and simplicial fan. 

Finally, we prove $\mathrm{PC}(\Sigma) =\Pi$. By (\ref{2.3.1}), every member $P$ of $\Pi$ does not generate a cone in $\Sigma$, but by (ii) in Definition 5 every proper subset of $P$ generates a cone in $\Sigma$. Therefore, $\Pi$ is a subset of $\mathrm{PC}(\Sigma)$. Conversely, let $S$ be in $\mathrm{PC}(\Sigma)$. Suppose that any member of $\Pi$ is not contained in $S$. Then $\langle S \rangle _{\ge 0} \in \Sigma$ by (\ref{2.3.1}), which contradicts that $S$ is a primitive collection of $\Sigma$. Hence there is a set $P \in \Pi$ such that $P \subset S$. If $P \subsetneq S$, then $\langle P \rangle _{\ge 0} \in \Sigma$ since $S$ is primitive (see Definition 4). However $P$ is primitive since we have already shown $\Pi \subset \mathrm{PC}(\Sigma)$. Thus $S$ must be equal to $P$ i.e., $S \in \Pi$. \ $\square$
\end{prf}

Here, we define the pseudo-fan corresponding to moduli space $\widetilde{Mp}_{0,2} (\mathbb{P}^1 \times \mathbb{P}^1 , (d_1 ,d_2))$ as follows:
\begin{defi}
\textit{\ Let $\nu ,\nu_1 ,\nu_2$ be members of $\mathrm{Col} (d_1 ,d_2 )$. Then we define
\begin{equation}
P_{\nu} := \begin{cases} 
\{ v_{a_0 ^1} ,v_{a_0 ^2} \} & (\nu = a_0) \\
\{ v_{a_{d_1} ^1} ,v_{a_{d_1} ^2} \} &(\nu = a_{d_1}) \\
\{ v_{b_0 ^1} ,v_{b_0 ^2} \} & (\nu = b_0) \\
\{ v_{b_{d_2} ^1} ,v_{b_{d_2} ^2} \} & (\nu = b_{d_2}) \\
\emptyset & (otherwise),  
              \end{cases}
\end{equation}
\begin{equation}
Q_{(\nu _1 ,\nu _2 )} := \begin{cases}
\{  v_{a_{i} ^1} ,v_{a_{i} ^2} ,v_{u_{(i,j)}} \} & ((\nu _1 ,\nu _2) = (a_i , u_{(i,j)}) \ (i=1 \dots ,d_1 -1 , j=0 \dots ,d_2)) \\
\{  v_{b_{j} ^1} ,v_{b_{j} ^2} ,v_{u_{(i,j)}} \} & ((\nu _1 ,\nu _2) = (b_j , u_{(i,j)}) \ (j=1 \dots ,d_2 -1 , i=0 \dots ,d_1)) \\
\{ v_{u_{(i,j)}} ,v_{u{(k,l)}} \} & ((\nu _1 ,\nu _2) = (u_{(i,j)} , u_{(k,l)}) \ (0 \le i < k \le d_1 ,0 \le l < j \le d_2)) \\
\emptyset & (otherwise), 
                               \end{cases}             
\end{equation}
and set
\begin{equation}
\Pi _{(d_1 ,d_2 )} :=\left( \{ P_{\nu} | \nu \in \mathrm{Col}(d_1 ,d_2) \} \cup \{ Q_{(\nu _1 ,\nu _2 )} | \nu_1 , \nu_2 \in \mathrm{Col} (d_1 ,d_2 ) \} \right) \setminus \{ \emptyset \},
\end{equation}
\begin{equation}
\Sigma _{(d_1 ,d_2)} := \Sigma(\mathrm{Ver}_{(d_1 ,d_2 )}, \Pi _{(d_1 ,d_2 )}).        
\end{equation}
}
\end{defi}
Then the triplet $(\Sigma _{(d_1 ,d_2)} , \mathrm{Ver}_{(d_1 ,d_2 )}, \Pi _{(d_1 ,d_2 )} )$ is clearly a pseudo-fan. By using Lemma 6, one of our main results in this paper can be stated as follows:
\begin{thm}
\textit{\ $(\Sigma _{(d_1 ,d_2)} , \mathrm{Ver}_{(d_1 ,d_2 )}, \Pi _{(d_1 ,d_2 )} )$ satisfies (MVC).}
\end{thm}
Instead of proving this theorem, we show the following lemma which is simpler than above theorem:
\begin{lem}
\textit{\ $(\Sigma _{(d_1 ,d_2 )} , \mathrm{Ver}_{(d_1 ,d_2 )}, \Pi _{(d_1 ,d_2 )})$ satisfies the following condition (AMVC) (Alternative Min--Value Condition):
\begin{description}
\item[(AMVC)] For any $y=(y_{\rho})_{\rho \in \mathrm{Ver} _{(d_1 ,d_2)}} \in \mathbb{R}^r$, there exists the unique $x=(x_{\mu})_{\mu \in \mathrm{Row}(d_1 ,d_2 )} \in \mathbb{R}^{r-n}$ such that
\begin{equation}
\mathrm{min} \left\{ -y_{\rho} + \sum _{\mu \in \mathrm{Row} (d_1 ,d_2 )} w_{(\mu ,\rho)} x_{\mu} \ \Bigg| \ \rho \in P \right\} =0 \ \ (\text{for any}\ P \in \Pi _{(d_1 ,d_2 )}). \label{AMVC}
\end{equation}
\end{description}
}
\end{lem}

Proof of this lemma is technical and complicated, and we will give the proof in Subsection 2.4. 
Assuming this lemma, we show Theorem 4 in the following. From Lemma 4, the map $V_{(d_1 ,d_2)} : \mathbb{Z}^r \rightarrow \mathbb{Z}^n$ is surjective and $V_{(d_1 ,d_2 )}{}^t\hspace{0.5pt} W_{(d_1 ,d_2)}={}^t\hspace{0.5pt}(W_{(d_1 ,d_2 )}{}^t\hspace{0.5pt} V_{(d_1 ,d_2)})=O$. Then for any $\beta \in \mathbb{R}^n$ and $x=(x_{\mu}) \in \mathbb{R}^{r-n}$, there exists $y =(y_{\rho}) \in \mathbb{R}^{r}$ such that 
\begin{align}
\beta &=V_{(d_1 ,d_2)} (-y + {}^t\hspace{0.5pt}W_{(d_1 ,d_2 )} x) \notag \\
        &= \sum _{\rho \in \mathrm{Ver}_{(d_1 ,d_2 )}} \left(-y_{\rho} + \sum _{\mu \in \mathrm{Row} (d_1 ,d_2 )} w_{(\mu ,\rho)} x_{\mu} \right) v_{\rho}.
\end{align}
Hence, we should put
\begin{equation}
\alpha _{\rho} = -y_{\rho} + \sum _{\mu \in \mathrm{Row} (d_1 ,d_2 )} w_{(\mu ,\rho)} x_{\mu},
\end{equation}
and by the condition (AMVC), we can take $x$ which satisfies 
\begin{equation}
\mathrm{min} \{ \alpha _{\rho} \ | \ \rho \in P \} = \mathrm{min} \left\{ -y_{\rho} + \sum _{\mu \in \mathrm{Row} (d_1 ,d_2 )} w_{(\mu ,\rho)} x_{\mu} \ \Bigg| \ \rho \in P \right\} =0
\end{equation}
for any $P \in \Pi _{(d_1 ,d_2 )}$. For the proof of uniqueness of $\alpha$ in (MVC), we assume that $\beta \in \mathbb{R}^n$ can be written as
\begin{equation}
\beta = \sum _{\rho} \alpha _{\rho} v_{\rho} =\sum _{\rho} \gamma _{\rho} v_{\rho}
\end{equation}
such that
\begin{equation}
\mathrm{min} \{ \alpha _{\rho} \ | \ \rho \in P \} = \mathrm{min} \{ \gamma _{\rho} \ | \ \rho \in P \} =0 \ \ (\text{for any}\ P \in \Pi _{(d_1 ,d_2 )})
\end{equation}
for some $\alpha =(\alpha _{\rho}), \gamma =(\gamma _{\rho}) \in \mathbb{R}^r$. By acting $\mathrm{Hom} _{\mathbb{Z}} (\bullet ,\mathbb{R})$ to (\ref{foundexactmat}), we obtain the following exact sequence:
\begin{equation}
\xymatrix@C=48pt{
0 \ar[r] & \mathbb{R}^{r-n} \ar[r]^-{{}^t\hspace{0.5pt} W_{(d_1 ,d_2)}} & \mathbb{R}^{r} \ar[r]^-{V_{(d_1 ,d_2)}} & \mathbb{R}^{n} \ar[r] & 0
}. \label{dexact}
\end{equation}
From the equality $V_{(d_1 ,d_2)} (\alpha -\gamma) =\sum _{\rho} (\alpha _{\rho} -\gamma _{\rho}) v_{\rho} =0$ and the above sequence, there exists $x=(x_{\mu}) \in \mathbb{R}^{r-n}$ such that $\alpha -\gamma ={}^t\hspace{0.5pt} W_{(d_1 ,d_2)} x$, i.e.,
\begin{equation}
\alpha _{\rho} =-(-\gamma _{\rho}) + \sum_{\mu} w_{(\mu ,\rho)} x_{\mu} \ (\text{for any $\rho \in \mathrm{Ver} _{(d_1 ,d_2 )}$}).
\end{equation}
Then $\mathrm{min} \{ -(-\gamma _{\rho}) + \sum_{\mu} w_{(\mu ,\rho)} x_{\mu} \ | \ \rho \in P \} =\mathrm{min} \{ \alpha _{\rho} \ | \ \rho \in P \} =0$. On the other hand, 
$\mathrm{min} \{ -(-\gamma _{\rho}) + \sum_{\mu} w_{(\mu ,\rho)} \cdot 0 \ | \ \rho \in P \} = \mathrm{min} \{ \gamma _{\rho} \ | \ \rho \in P \} =0$. Hence $x=0$ by uniqueness of $x$ in (AMVC) and therefore $\alpha = \gamma$.

\begin{ex}

\ The weight matrix $W_{(1,1)}$ is

\begin{equation}
W_{(1,1)} =(v_1 , \dots , v_{10}) = \left( \begin{array}{cccccccccc}
1 & 0 & 1 & 0 & 0 & 0 & 0 & 0 & 0 & -1 \\
0 & 1 & 0 & 1 & 0 & 0 & 0 & 0 & -1 & 0 \\ 
0 & 0 & 0 & 0 & 1 & 0 & 1 & 0 & -1 & 0 \\
0 & 0 & 0 & 0 & 0 & 1 & 0 & 1 & 0 & -1 \\
0 & 0 & 0 & 0 & 0 & 0 & 0 & 0 & 1 & 1 \\
\end{array}
\right)
\end{equation}
and the set $\Pi _{(1,1)}$ is 
\begin{equation}
\Pi _{(1,1)} =\{ \{ v_1 ,v_3 \} , \{ v_2 ,v_4 \} , \{ v_5 ,v_7 \} , \{ v_6 ,v_8 \} , \{ v_9 ,v_{10} \} \}.
\end{equation}
Thus, the system (\ref{AMVC}) is                                                                                                                                                                             
\begin{equation}
\begin{cases}
\mathrm{min} \{ z_1 ,z_3 \}=0 \\
\mathrm{min} \{ z_2 ,z_4 \}=0 \\
\mathrm{min} \{ z_5 ,z_7 \}=0 \\ 
\mathrm{min} \{ z_6 ,z_8 \}=0 \\
\mathrm{min} \{ z_9 ,z_{10} \}=0,
\end{cases}
\end{equation}
where
\begin{equation}
\left( \begin{array}{c}
z_1 \\
z_2 \\
z_3 \\
z_4 \\
z_5 \\
z_6 \\
z_7 \\
z_8 \\
z_9 \\
z_{10}
\end{array}
\right) :=
\left( \begin{array}{c}
-y_1 +x_1 \\
-y_2 +x_2 \\
-y_3 +x_1 \\
-y_4 +x_2 \\
-y_5 +x_3 \\
-y_6 +x_4 \\
-y_7 +x_3 \\
-y_8 +x_4 \\
-y_9 -x_2 -x_3 +x_5 \\
-y_{10} -x_1 -x_4 +x_5 
\end{array}
\right).
\end{equation}
This system can be solved easily;
\begin{equation}
\begin{cases}
x_1 =\mathrm{max} \{ y_1 , y_3 \} \\
x_2 =\mathrm{max} \{ y_2 , y_4 \} \\
x_3 =\mathrm{max} \{ y_5 , y_7 \} \\
x_4 =\mathrm{max} \{ y_6 , y_8 \} \\
x_5 =\mathrm{max} \{ y_9 +\mathrm{max} \{ y_2 , y_4 \} +\mathrm{max} \{ y_5 , y_7 \} , y_{10} +\mathrm{max} \{ y_1 , y_3 \} +\mathrm{max} \{ y_6 , y_8 \} \},
\end{cases}
\end{equation}
and therefore $\Sigma _{(1,1)}$ is a complete simplicial fan with $\mathrm{PC}(\Sigma _{(1,1)}) = \Pi _{(1,1)}$.
\end{ex}

\subsection{Quotient Construction}

In this subsection, we give the quotient construction of $X_{(d_1 ,d_2)} :=X({\Sigma _{(d_1 ,d_2)}})$. The exact sequence (\ref{1.2.8}) in case of $\Sigma =\Sigma _{(d_1 ,d_2)}$                                                                                             

\begin{equation}
\xymatrix@C=40pt{
0 \ar[r] & M \ar[r]^-{\varphi} & \mathbb{Z} ^{\Sigma _{(d_1 ,d_2)} (1)} \ar[r]^-{\psi} & A_{n-1} (X_{(d_1 ,d_2)}) \ar[r] & 0 
} \label{foundexact}
\end{equation}
gives us the way of calculation of Chow group $A_{n-1} (X_{(d_1 ,d_2 )})$.

\begin{lem}
\begin{equation}
A_{n-1} (X_{(d_1 ,d_2)}) \simeq \mathbb{Z}^{r-n}.
\end{equation}
\textit{In particular, the group  
$G= \mathrm{Hom} _{\mathbb{Z}} (A_{n-1} (X_{(d_1 ,d_2)}),\mathbb{C}^{\times})$
is isomorphic to
$(\mathbb{C}^{\times} )^{r-n}$.}
\end{lem}

\begin{prf}
\ Using the exact sequence (\ref{foundexact}),                                                                                                                                                    
\begin{align}
A_{n-1} (X_{(d_1 ,d_2)}) & \simeq \mathbb{Z} ^{\Sigma _{(d_1 ,d_2)} (1)} \big/ \mathrm{Ker} \, \psi \notag \\
                             & = \mathbb{Z} ^{\Sigma _{(d_1 ,d_2)} (1)} \big/ \mathrm{Im} \, \varphi \notag \\
                             & \simeq \mathbb{Z} ^{r} \Big/ \Bigl\langle (\langle e_i ,v_{j} \rangle)_{j=1 \dots ,r} \Big| i=1, \dots ,n \Bigr\rangle \notag \\
                             & \simeq \bigoplus _{j=1} ^{r} \mathbb{Z} x_j \Bigg/ \Biggl\langle \sum _{j=1} ^{r} \langle e_i ,v_{j} \rangle x_j \bigg| i=1, \dots ,n \Biggr\rangle 
\end{align}
($x_1 ,\dots ,x_r$ are formal symbols) and in this group,
\begin{equation}
\sum _{j=1} ^{r} \langle e_i ,v_{j} \rangle x_j =0 \ (i=1, \dots ,n) \Longleftrightarrow V_{(d_1 ,d_2)} 
\begin{pmatrix}
x_1 \\
x_2 \\
\vdots \\
x_r
\end{pmatrix}
=O.
\end{equation} 

From this fact and definition of $V_{(d_1 ,d_2)}$, we can easily see that $A_{n-1} (X_{(d_1 ,d_2)})$ is generated by $r-n$ elements. $\square$
\end{prf}

Let us denote the $T_{N}$--invariant divisor that corresponds to $v _j \in \mathrm{Ver} _{(d_1 ,d_2)}$ by $D_j$.

\begin{lem}
\textit{\ We can choose $\mathbb{Z}$--bases $E_1 , \dots ,E_{r-n}$ of $A_{n-1} (X_{(d_1 ,d_2 )})$ such that  
\begin{equation}
[D_j]=\sum _{i=1}^{r-n} w_{(i,j)} E_i \quad (j=1 ,\dots ,r). 
\end{equation}
}
\end{lem}

\begin{prf}       
\ By Lemma 4, there exists $\beta _1 , \dots ,\beta _{r-n} \in \mathbb{Z}^r$ such that $W_{(d_1 ,d_2)} \beta _i = e_i^{(r-n)} \ (i=1, \dots ,r-n)$, where $e_i^{(l)}$ is a $i$--th canonical base of $\mathbb{Z}^l$. We set $E_i := \psi (\beta _i)$ \ ($\psi$ is the map in exact sequence (\ref{foundexact})). Then                                                                             
\begin{align}
W_{(d_1 ,d_2 )} e_j ^{(r)} &= \sum_{i=1}^{r-n} w_{(i,j)} e_{i}^{(r-n)} \notag \\ 
             &= \sum_{i=1}^{r-n} w_{(i,j)} W_{(d_1 ,d_2 )} \beta_{i} \notag \\
             &=W_{(d_1 ,d_2 )} \Biggl(\sum_{i=1}^{r-n} w_{(i,j)} \beta_{i} \Biggr), \notag \\
\end{align}
and by Lemma 4 and (\ref{foundexact}), $e_{i}^{(r)} -\sum_{i=1}^{r-n} w_{(i,j)} \beta_{i} \in \mathrm{Ker} \, W_{(d_1 ,d_2 )} =\mathrm{Ker} \, \psi$. Thus
\begin{align}
[D_j] &=\psi (e_{j}^{(r)}) \notag \\
       &=\psi \left( \sum_{i=1}^{r-n} w_{(i,j)} \beta_{i} \right) \notag \\
       &=\sum_{i=1}^{r-n} w_{(i,j)} \psi (\beta_i) \notag \\
       &=\sum_{i=1}^{r-n} w_{(i,j)} E_i . \label{lemfo}
\end{align}

In order to complete the proof, we must check that $E_1 ,\dots ,E_{r-n}$ are $\mathbb{Z}$--bases of $A_{n-1} (X_{(d_1 ,d_2)})$. It is clear that $E_1 ,\dots ,E_{r-n}$ generate $A_{n-1} (X_{(d_1 ,d_2)})$ by Lemma 4 and (\ref{lemfo}). If $\sum _{i=1}^{r-n} c_{i} E_i =0 \ (c_{i} \in \mathbb{Z})$, we have, 
\[ \psi \Biggl( \sum_{i=1}^{r-n} c_{i} \beta_{i} \Biggr) =\sum_{i=1}^{r-n} c_{i} E_{i} =0. \]
Since $\mathrm{Ker} \, W_{(d_1 ,d_2 )} =\mathrm{Ker} \, \psi$, we also have, 
\[ \sum_{i=1}^{r-n} c_{i} e_{i}^{(r-n)} = W_{(d_1 ,d_2 )} \Biggl( \sum_{i=1}^{r-n} c_{i} \beta_{i} \Biggr) =0. \]
Then we can conclude $c_1 = \dots =c_{r-n} =0.$ $\Box$
\end{prf}

In our situation, $G$ acts on $\mathbb{C}^{r}$ as follows (see (\ref{1.2.11}) in Subsection 1.2):
\begin{equation}
G \times \mathbb{C}^{r} \ni (g,(x_{j})_{j}) \longmapsto (g([D_{j}])x_{j})_{j} \in \mathbb{C}^{r}.
\end{equation}
From Lemma 9, we obtain,
\begin{align}
g([D_j]) &=g \Biggl( \sum_{i=1}^{r-n} w_{(i,j)} E_{i} \Biggr) \notag \\
          &=\prod_{i=1}^{r-n} g(E_{i})^{w_{(i,j)}},
\end{align}
and since $G$ is identified with $(\mathbb{C}^{\times})^{r-n}$ by the correspondence: $g \longleftrightarrow (\lambda _1, \dots ,\lambda_{r-n} ) = (g(E_1) ,\dots ,g(E_{r-n}))$, we reach the following theorem by using Lemma 3 and Theorem 2:                                                                                                                                                            
\begin{thm}
\[ X_{(d_1 ,d_2)}=(\mathbb{C}^{\Sigma _{(d_1 ,d_2 )} (1)} - Z(\Sigma _{(d_1 ,d_2) }))/G = \widetilde{Mp}_{0,2} (\mathbb{P}^1 \times \mathbb{P}^1 , (d_1 ,d_2)). \]
\end{thm}
Combining  the above theorem with  Theorem 1, we obtain the main result of this paper.                                                                                                               
\begin{thm}
\textit{\ The moduli space $\widetilde{Mp}_{0,2} (\mathbb{P}^1 \times \mathbb{P}^1 , (d_1 ,d_2))$ is a complete and simplicial toric variety. Especially, it is compact orbifold.}
\end{thm}

\subsection{Proof of Lemma 7}                                                                                                                                                  

We denote the set $( \{ 0 ,\dots , d_1 \} \times \{ 0,\dots ,d_2 \} ) \setminus \{ (0,0),(d_1 ,d_2 ) \}$ by $I$ and define
\begin{equation}
J_{(i,j)}^{+} :=\{ (k,l) \in \{  0, \dots ,d_1 \} \times \{ 0,\dots ,d_2 \} \ \big| \ 0 \le i < k \le d_1 , 0 \le l < j \le d_2 \},
\end{equation}
\begin{equation}
J_{(i,j)}^{-} :=\{ (k,l) \in \{  0, \dots ,d_1 \} \times \{ 0,\dots ,d_2 \} \ \big| \ 0 \le k < i \le d_1 , 0 \le j < l \le d_2 \},
\end{equation}
and
\begin{equation}
J_{(i,j)} :=J_{(i,j)}^{+} \cup J_{(i,j)}^{-}
\end{equation}
for $(i,j) \in I$.

Before solving the system (\ref{AMVC}), we show the following simple lemma.                                                                                                                                  

\begin{lem}
\textit{\ Let $a ,b_1 ,\dots , b_m$ be real numbers. Then $\mathrm{min} \{ a ,b_i \} =0$ for all $i=1,\dots ,m$ if and only if $\mathrm{min} \{ a , b_1 + \dots + b_m \} =0$ and $b_1 ,\dots ,b_m \ge 0$.
}
\end{lem}

\begin{prf}
\ If $a=0$, this lemma is obvious. Hence we assume that $a \neq 0$. If $\mathrm{min} \{ a ,b_i \} =0$ for all $i=1,\dots ,m$, then $a>0$ and $b_1 ,\dots ,b_m =0$ i.e., $b_1 + \dots  + b_m =0$. Conversely, if $\mathrm{min} \{ a , b_1 + \dots + b_m \} =0$ and $b_1 ,\dots ,b_m \ge 0$, then $a>0$ and $b_1 + \dots + b_m =0$. Since $b_1 ,\dots ,b_m \ge 0$, $b_1 \dots ,b_m$ must be equal to $0$. \ $\square$ 
\end{prf}

Now, we solve the system (\ref{AMVC}) using this lemma and Lemma 5. In case of $P_{a_0} ,P_{a_{d_1}} ,P_{b_0} ,P_{b_{d_2}}$, the crresponding min--value equation are          
\begin{equation}
\begin{cases}
\mathrm{min} \{ -y_{a_0 ^1} +x_{z_0} , -y_{a_0 ^2} +x_{z_0} \} =0 \ \ (P_{a_0}) \\
\mathrm{min} \{ -y_{a_{d_1} ^1} +x_{z_{d_1}} , -y_{a_{d_1} ^2} +x_{z_{d_1}} \} =0 \ \ (P_{a_{d_1}}) \\ 
\mathrm{min} \{ -y_{b_0 ^1} +x_{w_0} , -y_{b_0 ^2} +x_{w_0} \} =0 \ \ (P_{b_0}) \\
\mathrm{min} \{ -y_{b_{d_2} ^1} +x_{w_{d_2}} , -y_{b_{d_2} ^2} +x_{w_{d_2}} \} =0 \ \ (P_{b_{d_2}}),
\end{cases} \label{2.4.4}
\end{equation}
and these are easily solved as follows:
\begin{equation}
\begin{cases}
x_{z_0} =\mathrm{max} \{ y_{a_0 ^1} ,y_{a_0 ^2} \} \\
x_{z_{d_1}} =\mathrm{max} \{ y_{a_{d_1} ^1} ,y_{a_{d_1} ^2} \} \\
x_{w_0} =\mathrm{max} \{ y_{b_0 ^1} ,y_{b_0 ^2} \} \\
x_{w_{d_2}} =\mathrm{max} \{ y_{b_{d_2} ^1} ,y_{b_{d_2} ^2} \}.
\end{cases} \label{2.4.5}
\end{equation}
For each $i=1 ,\dots ,d_1 -1$, equations corresponding to $Q_{(a_i ,u_{(i,j)})}$ ($j=0 ,\dots ,d_2$) are combined into,
\begin{equation}
\begin{cases}
\mathrm{min} \{ -y_{a_i}+x_{z_i} , -y_{(i,0)} - \dots -y_{(i,d_2 )} -x_{z_{i-1}} +2x_{z_i} -x_{z_{i+1}} \} =0, \\
-y_{(i,j)} + \sum_{\mu \in \mathrm{Row} (d_1 ,d_2)} w_{(\mu ,u_{(i,j)})} x_{\mu} \ge 0 \ \ (j=0, \dots ,d_2),
\end{cases} \label{2.4.6}
\end{equation} 
by using (\ref{2.1.24}) in Lemma 5, Lemma 10, and the equalities $\mathrm{min} \{ a,b,c \} =\mathrm{min} \{ \mathrm{min} \{ a,b \} ,c \}$, $\mathrm{min} \{ -a+x ,-b+x \} =x-\mathrm{max} \{ a,b \}$. Note here that we use the abbreviation: $y_{(i,j)} :=y_{u_{(i,j)}}$, $y_{a_i} := \mathrm{max} \{ y_{a_i ^1} ,y_{a_i ^2} \}$.                                                                                                                  
Since the inequalities in the second line of (\ref{2.4.6}) appear later again (see (\ref{2.4.8}) and the second line of (\ref{2.4.11}) and (\ref{2.4.12})), we only need to consider the following system which consists of the min--value equations presented in the first lines of (\ref{2.4.6}).                                                           
\begin{equation}
\begin{cases}
\mathrm{min} \{ -y_{a_1}+x_{z_1} , -y_{(1,0)} - \dots -y_{(1,d_2 )} -x_{z_{0}} +2x_{z_1} -x_{z_{2}} \} =0 \\
\mathrm{min} \{ -y_{a_2}+x_{z_2} , -y_{(2,0)} - \dots -y_{(2,d_2 )} -x_{z_{1}} +2x_{z_2} -x_{z_{3}} \} =0 \\
\dotsb \\
\mathrm{min} \{ -y_{a_{d_1 -1}}+x_{z_{d_1 -1}} , -y_{(d_1 -1,0)} - \dots -y_{(d_1 -1,d_2 )} -x_{z_{d_1 -2}} +2x_{z_{d_1 -1}} -x_{z_{d_1}} \} =0,
\end{cases} \label{2.4.7}
\end{equation}
($x_{z_0} :=\mathrm{max} \{ y_{a_0 ^1} ,y_{a_0 ^2} \}$, $x_{z_{d_1}} :=\mathrm{max} \{ y_{a_{d_1} ^1} ,y_{a_{d_1} ^2} \}$). 
It is already shown by Saito \cite{S} that this system has unique solution of $x_{z_0} ,\dots ,x_{z_{d_1}}$ for fixed $y$'s. Similarly, we can also find unique $x_{w_{1}} , \dots ,x_{w_{d_2 -1}}$ solutions from the system corresponding to $Q_{(b_j ,u_{(i,j)})}$ ($j=1, \dots ,d_2 -1 ,i=0,\dots ,d_1$). At this stage, we denote these solutions by $x_{z_{i}} ^{\prime}$, $x_{w_{j}} ^{\prime}$ and define,
\begin{equation}
y_{(i,j)} ^{\prime} :=y_{(i,j)} -\sum _{p=0} ^{d_1} w_{(z_{p} , u_{(i,j)})} x_{z_{p}} ^{\prime} -\sum _{q=0} ^{d_2} w_{(w_{q} , u_{(i,j)})} x_{w_{q}} ^{\prime} \ \ ((i,j) \in I). \label{2.4.8}
\end{equation}
The second line of  (\ref{2.4.6}) tells us that                                                                                                                                                                                 
\begin{equation}
\sum _{j=0} ^{d_2} y_{(i,j)} ^{\prime} \le 0 \ \ (i=1 ,\dots ,d_1). \label{2.4.9}
\end{equation}
The inequalities corresponding to $Q_{(b_j ,u_{(i,j)})}$ ($j=1, \dots ,d_2 -1 ,i=0,\dots ,d_1$) also gives,
\begin{equation}
\sum _{i=0} ^{d_1} y_{(i,j)} ^{\prime} \le 0 \ \ (j=1 ,\dots ,d_2). \label{2.4.10}
\end{equation}

Next, we find the solution $x_{(i,j)} :=x_{f_{(i,j)}}$. Let $(i,j)$ be fixed member of $I$. From (\ref{2.1.27}) in Lemma 5 and Lemma 10, equations corresponding to $Q_{(u_{(i,j)} ,u_{(k,l)})}$ ($(k,l) \in J_{(i,j)}^{+}$) are combined into,                                                                                                                                                                                             
\begin{equation}
\begin{cases}
\mathrm{min} \left\{ -y_{(i,j)}^{\prime} +\sum _{p=1} ^{d_1} \sum _{q=1} ^{d_2} w_{(f_{(p,q)} ,u_{(i,j)})} x_{(p,q)} , -\sum_{p=i+1}^{d_1} \sum_{q=0}^{j-1} y_{(p,q)}^{\prime} +x_{(i+1,j)} \right\} =0, \\
-y_{(k,l)}^{\prime} +\sum _{p=1} ^{d_1} \sum _{q=1} ^{d_2} w_{(f_{(p,q)} ,u_{(k,l)})} x_{(p,q)} \ge 0 \ \ (\text{for all}\ (k,l) \in J_{(i,j)}^{+}),
\end{cases} \label{2.4.11}
\end{equation}
and the ones corresponding to $Q_{(u_{(i,j)} ,u_{(k,l)})}$ ($(k,l) \in J_{(i,j)}^{-}$) are combined into,
\begin{equation}
\begin{cases}
\mathrm{min} \left\{ -y_{(i,j)}^{\prime} +\sum _{p=1} ^{d_1} \sum _{q=1} ^{d_2} w_{(f_{(p,q)} ,u_{(i,j)})} x_{(p,q)} , -\sum_{p=0}^{i-1} \sum_{q=j+1}^{d_1} y_{(p,q)}^{\prime} +x_{(i,j+1)} \right\} =0, \\
-y_{(k,l)}^{\prime} +\sum _{p=1} ^{d_1} \sum _{q=1} ^{d_2} w_{(f_{(p,q)} ,u_{(k,l)})} x_{(p,q)} \ge 0 \ \ (\text{for all}\ (k,l) \in J_{(i,j)}^{-}),
\end{cases} \label{2.4.12}
\end{equation}
where we put $x_{(s,t)}:=0$ ($(s,t) \notin \{ 1,\dots ,d_1 \} \times \{ 1,\dots ,d_2 \}$) and $y_{(p,q)}^{\prime} :=0$ ($(p,q) \notin I$). Since
\begin{equation}
w_{(f_{(p,q)} ,u_{(s,t)})}= \begin{cases}
                         -1 \ \ ((p,q)=(s,t),(s+1,t+1)) \\
                         1 \ \ ((p,q)=(s+1,t),(s,t+1)) \\
                         0 \ \ (otherwise)
                         \end{cases}
\end{equation}
for $p=1,\dots , d_1$,$q=1 \dots ,d_2$ and $(s,t) \in I$, (\ref{2.4.11}) and (\ref{2.4.12}) are further combined into,                                                                                                    
\begin{equation}
\begin{cases}
\mathrm{min} \{ -y_{(i,j)}^{\prime} -x_{(i,j)} +x_{(i+1,j)} +x_{(i,j+1)} -x_{(i+1,j+1)} , \\
\qquad -\Delta_{(i+1,0)}^{(d_1 ,j-1)} -\Delta_{(0,j+1)}^{(i-1,d_2 )} +x_{(i+1,j)}+x_{(i,j+1)} \} =0, \\
-y_{(k,l)}^{\prime} -x_{(k,l)} +x_{(k+1,l)} +x_{(k,l+1)} -x_{(k+1,l+1)} \ge 0 \ \ (\text{for all}\ (k,l) \in J_{(i,j)}),
\end{cases}
\end{equation}
where we introduce the following notation:
\begin{equation}
\Delta _{(p,q)}^{(s,t)} := \begin{cases}
                             \sum_{k=p}^{s} \sum_{l=q}^{t} y_{(p,q)}^{\prime} \ \ (p \le s , q \le t) \\
                             0 \ \ (otherwise).
                             \end{cases}
\end{equation}
Since the l.h.s. of the inequalities in the second line of (\ref{2.4.11}) or (\ref{2.4.12}) corresponding to $Q_{(u_{(i,j)},u_{(k,l)})}$ also appear  
in the min--value equations in the first line of (\ref{2.4.11}) or (\ref{2.4.12}) corresponding to another $Q_{(u_{(i,j)},u_{(k,l)})}$, 
we can omit the second lines of (\ref{2.4.11}) and (\ref{2.4.12}).
Therefore, the following system of min--value equations is equivalent to equations corresponding to $Q_{(u_{(i,j)},u_{(k,l)})} \ ((k,l) \in J_{(i,j)})$ in (\ref{AMVC}):                                                                              
\begin{multline}
\mathrm{min} \{ -y_{(i,j)}^{\prime} -x_{(i,j)} +x_{(i+1,j)} +x_{(i,j+1)} -x_{(i+1,j+1)} , \\
-\Delta_{(i+1,0)}^{(d_1 ,j-1)} -\Delta_{(0,j+1)}^{(i-1,d_2 )} +x_{(i+1,j)}+x_{(i,j+1)} \} =0 \ \ ((i,j) \in I).
\end{multline}
The above system is further rewritten as follows:
\begin{equation}
x_{(i+1,j)} +x_{(i,j+1)}=\mathrm{max} \{ y_{(i,j)}^{\prime} +x_{(i,j)}+x_{(i+1,j+1)} , \Delta_{(i+1,0)}^{(d_1 ,j-1)} +\Delta_{(0,j+1)}^{(i-1,d_2 )} \}
\ \ ((i,j) \in I). \label{2.4.17}
\end{equation}
\begin{lem}
\textit{\ The unique solution of (\ref{2.4.17}) is given by,                                                                                                                                                      
\begin{equation}
x_{(i,j)} =\mathrm{max} \{ \Delta_{(i,0)}^{(d_1 ,j-1)},\Delta_{(0,j)}^{(i-1,d_2)} \} \ \ ((i,j) \in I). \label{2.4.18}
\end{equation}
}
\end{lem}

\begin{prf}
\ First we check that (\ref{2.4.18}) satisfies the equation (\ref{2.4.17}). Since $\mathrm{max} \{ a,b \} + \mathrm{max} \{ c,d \} =\mathrm{max} \{ a+c,a+d,b+c,b+d \}$, we obtain,
\begin{align}
x_{(i+1,j)} + x_{(i,j+1)} &= \mathrm{max} \{ \Delta_{(i+1,0)}^{(d_1 ,j-1)},\Delta_{(0,j)}^{(i,d_2)} \} + \mathrm{max} \{ \Delta_{(i,0)}^{(d_1 ,j)},\Delta_{(0,j+1)}^{(i-1,d_2)} \} \notag \\
                          &= \mathrm{max} \{ \Delta_{(i+1,0)}^{(d_1 ,j-1)} +\Delta_{(i,0)}^{(d_1 ,j)},\Delta_{(i+1,0)}^{(d_1 ,j-1)} + \Delta_{(0,j+1)}^{(i-1,d_2)}, \notag \\
                          &\ \quad \Delta_{(0,j)}^{(i,d_2)} +\Delta_{(i,0)}^{(d_1 ,j)} , \Delta_{(0,j)}^{(i,d_2)} +\Delta_{(0,j+1)}^{(i-1,d_2)} \}.
\end{align}
At this stage, note that
\begin{align}
\Delta_{(0,j)}^{(i,d_2)} +\Delta_{(i,0)}^{(d_1 ,j)} &=\Delta_{(i+1,0)}^{(d_1 ,j-1)} + \Delta_{(0,j+1)}^{(i-1,d_2)} + \sum _{j=0} ^{d_2} y_{(i,j)} ^{\prime} + \sum _{i=0} ^{d_1} y_{(i,j)} ^{\prime} \notag \\
                                                       &\le \Delta_{(i+1,0)}^{(d_1 ,j-1)} + \Delta_{(0,j+1)}^{(i-1,d_2)},
\end{align}
holds by (\ref{2.4.9}) and (\ref{2.4.10}). Therefore, we have,
\begin{equation}
x_{(i+1,j)} +x_{(i,j+1)} =  \mathrm{max} \{ \Delta_{(i+1,0)}^{(d_1 ,j-1)} +\Delta_{(i,0)}^{(d_1 ,j)},\Delta_{(i+1,0)}^{(d_1 ,j-1)} + \Delta_{(0,j+1)}^{(i-1,d_2)}, \Delta_{(0,j)}^{(i,d_2)} +\Delta_{(0,j+1)}^{(i-1,d_2)} \}.
\end{equation}
On the other hand, the r.h.s. of (\ref{2.4.17}) is rewritten as follows.                                                                                                                                                                               
\begin{align}
& \mathrm{max} \{ y_{(i,j)}^{\prime} +x_{(i,j)}+x_{(i+1,j+1)} , \Delta_{(i+1,0)}^{(d_1 ,j-1)} +\Delta_{(0,j+1)}^{(i-1,d_2 )} \} \notag \\
&= \mathrm{max} \bigl\{ y_{(i,j)}^{\prime} + \mathrm{max} \{ \Delta_{(i,0)}^{(d_1 ,j-1)},\Delta_{(0,j)}^{(i-1,d_2)} \} + \mathrm{max} \{ \Delta_{(i+1,0)}^{(d_1 ,j)},\Delta_{(0,j+1)}^{(i,d_2)} \} , \notag \\
&\ \quad \Delta_{(i+1,0)}^{(d_1 ,j-1)} +\Delta_{(0,j+1)}^{(i-1,d_2 )} \bigr\} \notag \\
&= \mathrm{max} \bigl\{ y_{(i,j)}^{\prime} + \mathrm{max} \{ \Delta_{(i,0)}^{(d_1 ,j-1)} +\Delta_{(i+1,0)}^{(d_1 ,j)},\Delta_{(i,0)}^{(d_1 ,j-1)} +\Delta_{(0,j+1)}^{(i,d_2)} , \notag \\
&\ \quad \Delta_{(0,j)}^{(i-1,d_2)} + \Delta_{(i+1,0)}^{(d_1 ,j)} , \Delta_{(0,j)}^{(i-1,d_2)} + \Delta_{(0,j+1)}^{(i,d_2)} \} ,  \Delta_{(i+1,0)}^{(d_1 ,j-1)} +\Delta_{(0,j+1)}^{(i-1,d_2 )} \bigr\} \notag \\
&= \mathrm{max} \bigl\{ y_{(i,j)}^{\prime} +\Delta_{(i,0)}^{(d_1 ,j-1)} +\Delta_{(i+1,0)}^{(d_1 ,j)},y_{(i,j)}^{\prime} +\Delta_{(i,0)}^{(d_1 ,j-1)} +\Delta_{(0,j+1)}^{(i,d_2)} , \notag \\
&\ \quad y_{(i,j)}^{\prime} +\Delta_{(0,j)}^{(i-1,d_2)} + \Delta_{(i+1,0)}^{(d_1 ,j)} , y_{(i,j)}^{\prime} +\Delta_{(0,j)}^{(i-1,d_2)} + \Delta_{(0,j+1)}^{(i,d_2)} , \Delta_{(i+1,0)}^{(d_1 ,j-1)} +\Delta_{(0,j+1)}^{(i-1,d_2 )} \bigr\} \notag \\
&= \mathrm{max} \Biggl\{ \Delta_{(i,0)}^{(d_1 ,j)} +\Delta_{(i+1,0)}^{(d_1 ,j-1)},\Delta_{(i+1,0)}^{(d_1 ,j-1)} +\Delta_{(0,j+1)}^{(i-1,d_2)}+\sum _{j=0} ^{d_2} y_{(i,j)} ^{\prime} , \notag \\
&\ \quad \Delta_{(0,j+1)}^{(i-1,d_2)} + \Delta_{(i+1,0)}^{(d_1 ,j-1)} + \sum _{i=0} ^{d_1} y_{(i,j)} ^{\prime} ,\Delta_{(0,j)}^{(i,d_2)} + \Delta_{(0,j+1)}^{(i-1,d_2)} , \Delta_{(i+1,0)}^{(d_1 ,j-1)} +\Delta_{(0,j+1)}^{(i-1,d_2 )} \Biggr\} \notag \\
&= \mathrm{max} \bigl\{ \Delta_{(i,0)}^{(d_1 ,j)} +\Delta_{(i+1,0)}^{(d_1 ,j-1)},\Delta_{(0,j)}^{(i,d_2)} + \Delta_{(0,j+1)}^{(i-1,d_2)} , \Delta_{(i+1,0)}^{(d_1 ,j-1)} +\Delta_{(0,j+1)}^{(i-1,d_2 )} \bigr\}.
\end{align}
Hence $x_{(i,j)} =\mathrm{max} \{ \Delta_{(i,0)}^{(d_1 ,j-1)},\Delta_{(0,j)}^{(i-1,d_2)} \}$ satisfy the system (\ref{2.4.17}). Uniqueness
of the solution is obvious because $x_{(1,d_2)}$ is uniquely determined by the equation in $(i,j)=(0,d_2)$: 
\begin{equation}
x_{(1,d_2)} =\mathrm{max} \{ y_{(0,d_2)}^{\prime} , \Delta_{(1,0)}^{(d_1 ,d_2 -1)} \}
\end{equation}
and the remaining $x_{(i,j)}$'s are inductively determined by other equations in (\ref{2.4.17}). \ $\square$
\end{prf}
In this way, we have constructed unique solution of (\ref{AMVC}) for arbitrary $y=(y_{\rho})_{\rho \in \mathrm{Ver} _{(d_1 ,d_2)}} \in \mathbb{R}^r$ and have completed the proof of Lemma 7.

\section{Chow Ring and Poincar\'e Polynomial of $\widetilde{Mp}_{0,2} (\mathbb{P}^1 \times \mathbb{P}^1 , (d_1 ,d_2))$}

In this section, we compute Chow ring and Poincar\'e polynomial of $\widetilde{Mp}_{0,2} (\mathbb{P}^1 \times \mathbb{P}^1 , (d_1 ,d_2))$ using Theorem 3. For the sake of simplicity,  we use the order of elements of $\mathrm{Ver} _{(d_1 ,d_2 )}(1)$
presented in (\ref{vertex}), i.e., we denote the $i$-th column vector of $V_{(d_1 ,d_2 )}$ by $v_i$.
\begin{thm}
\begin{align}
A_{\mathbb{Q}} ^{\ast} (X_{(d_1 ,d_2)}) & \simeq H^{\ast} (X_{(d_1 ,d_2)} ,\mathbb{Q}) \notag \\
                                                & \simeq \mathbb{Q} [h_1 ,\dots ,h_{r-n}] \Bigg/ \Biggl\langle \prod _{v_j \in S} (w_{(1,j)} h_1 + \dots + w_{(r-n,j)} h_{r-n}) \bigg| S \in \Pi _{(d_1 ,d_2 )} \Biggr\rangle . \label{3.1}
\end{align}
\end{thm}

\begin{prf}
\ First we prove that $\mathbb{Q}[x_1 , \dots ,x_r] / I(\Sigma _{(d_1 ,d_2 )}) \simeq \mathbb{Q} [h_1 ,\dots ,h_{r-n}]$ with identification
\begin{equation}
x_{i} = \sum_{k=1}^{r-n} w_{(k,i)} h_{k}.  \label{3.2}
\end{equation}
for $i=1, \dots ,r$. We define a $\mathbb{Q}$-algebra homomorphism $\theta$ from $\mathbb{Q}[x_1 , \dots ,x_r]$ to $\mathbb{Q} [h_1 ,\dots ,h_{r-n}]$ by $x_{i} \longmapsto \sum_{k=1}^{r-n} w_{(k,i)} h_{k}$. Since $V_{(d_1 ,d_2 )} {}^t\hspace{0.5pt} W_{(d_1 ,d_2)}=O$, $\theta$ induces a $\mathbb{Q}$-algebra homomorphism $\overline{\theta}$ from quotient ring $\mathbb{Q}[x_1 , \dots , x_r]/I(\Sigma _{(d_1 ,d_2 )})$ to $\mathbb{Q} [h_1 ,\dots ,h_{r-n}]$. By using the exact sequence (\ref{dexact}), we can obtain the following exact sequence, 
\begin{equation}
\xymatrix@C=48pt{
0 \ar[r] & \mathbb{Q}^{r-n} \ar[r]^-{{}^t\hspace{0.5pt} W_{(d_1 ,d_2)}} & \mathbb{Q}^{r} \ar[r]^-{V_{(d_1 ,d_2)}} & \mathbb{Q}^{n} \ar[r] & 0
.}
\end{equation}
Obviously, this sequence splits (i.e., $\mathbb{Q}^r \simeq \mathbb{Q}^{r-n} \oplus \mathbb{Q}^{n}$). Then we can take two matrices $A=(a_{ki}) \in M(r-n,r;\mathbb{Q}),B=(b_{ij}) \in M(r,n;\mathbb{Q})$ such that ${}^t\hspace{0.5pt} W_{(d_1 ,d_2)} A + B V_{(d_1 ,d_2 )} =I_{r}$, $A {}^t\hspace{0.5pt} W_{(d_1 ,d_2)}=I_{r-n}$, and $V_{(d_1 ,d_2 )}B=I_{n}$. Using the matrix $A$, we can define a $\mathbb{Q}$-algebra homomorphism $\tau$ from $\mathbb{Q} [h_1 ,\dots ,h_{r-n}]$ to $\mathbb{Q}[x_1 , \dots ,x_r]/I(\Sigma _{(d_1 ,d_2 )})$ by $h_k \longmapsto [\sum_{i=1}^{r} a_{ki} x_{i}]$. Then 
\begin{align}
\tau \circ \overline{\theta} ([x_{i}]) &= \sum_{k=1}^{r-n} w_{(k,i)} \sum_{l=1}^{r} a_{kl} [x_{l}] \notag \\
                                                  &=\sum_{l=1}^{r} ({}^t\hspace{0.5pt} W_{(d_1 ,d_2)} A)_{il} [x_{l}] \notag \\
                                                  &=\sum_{l=1}^{r} (I_{r} -B V_{(d_1 ,d_2 )})_{il} [x_{l}] \notag \\
                                                  &=[x_{i}] -\sum _{j=1} ^{n} b_{ij} \left[ \sum_{l=1}^{r} v_{jl} x_{l} \right] \notag \\
                                                  &=[x_{i}] \quad (\text{by the definition of $I(\Sigma)$}),
\end{align}
where $(C)_{ij}$ is $(i,j)$--component of matrix $C$. We can also show $\overline{\theta} \circ \tau (h_k) =h_k$ by using the relation $A {}^t\hspace{0.5pt} W_{(d_1 ,d_2)}=I_{r-n}$. Hence $\overline{\theta}$ is an isomorphism. Then the identification (\ref{3.2}) and the definition of Stanley--Reisner ideal directly lead us to (\ref{3.1}). \ $\square$
\end{prf}

\begin{ex}[$(d_1 ,d_2)=(1,1)$]
\[ A_{\mathbb{Q}} ^{\ast} (X_{(1,1)}) \simeq \mathbb{Q} [h_1 ,\dots , h_5] / \langle h_1 ^2 , h_2 ^2 , h_3 ^2 , h_4 ^2 , (-h_2 -h_3 +h_5 )(-h_1 -h_4 +h_5 ) \rangle \]
\end{ex}

\begin{ex}[$(d_1 ,d_2)=(2,1)$]
\begin{align}
A_{\mathbb{Q}} ^{\ast} (X_{(2,1)}) \simeq & \mathbb{Q} [h_1 ,\dots , h_7] / \langle h_1 ^2 , h_3 ^2 , h_4 ^2 , h_5 ^2 , h_2 ^2 (-h_1 +h_2 +h_6 -h_7) , h_2 ^2 (h_2 -h_3 -h_6 +h_7), \notag \\
& (-h_2 -h_4 +h_6)(-h_1 +h_2 +h_6 -h_7) ,(-h_2 -h_4 +h_6)(h_2 -h_3 -h_6 +h_7), \notag \\
& (h_2 -h_3 -h_6 +h_7)(-h_2 -h_5 +h_7) \rangle \notag
\end{align}
\end{ex}

Next, we compute Poincar\'e Polynomial of $\widetilde{Mp}_{0,2} (\mathbb{P}^1 \times \mathbb{P}^1 , (d_1 ,d_2))$ in case of $(d_1 ,d_2 )=(1,1),(2,1)$ by using Theorem 3. At first sight, it seems that we have to count the number $|\Sigma _{(d_1 ,d_2 )}(k)|$ for  $k=0 ,\dots, n$. But actually, thanks to the facts (iii) and (iv) in Theorem 3, we only have to count the numbers $|\Sigma _{(d_1 ,d_2 )}(0)|, \dots ,|\Sigma _{(d_1 ,d_2 )}((n-1)/2)| = |\Sigma _{(d_1 ,d_2 )}(d_1 +d_2 )|$. Moreover, $|\Sigma _{(d_1 ,d_2 )}(k)|$ equals to the number: 
\begin{equation}
\big| \{ S \ \big| \ S \subset \mathrm{Ver} _{(d_1 ,d_2 )}, |S|=k \text{, and all $P \in \Pi _{(d_1 ,d_2 )}$ are not contained in $S$} \} \big|.
\end{equation}
by definition of $\Sigma _{(d_1 ,d_2 )}$ (see (\ref{2.3.1})).

\begin{description}
\item[(i) \ ($(d_1 ,d_2 )=(1,1)$)]
Since $|\Sigma _{(1,1)} (0)|=1, |\Sigma _{(1,1)} (1)|=10, |\Sigma _{(1,1)} (2)|=40$, Betti numbers are computed using Theorem 3 as follows:
\begin{align}
&b_{0} (X_{(1,1)}) = b_{10} (X_{(1,1)}) =|\Sigma _{(1,1)} (0)|=1, \notag \\
&b_{2} (X_{(1,1)}) = b_{8} (X_{(1,1)}) =|\Sigma _{(1,1)} (1)|-5|\Sigma _{(1,1)} (0)|=5, \notag \\
&b_{4} (X_{(1,1)}) = b_{6} (X_{(1,1)}) =|\Sigma _{(1,1)} (2)|-4|\Sigma _{(1,1)} (1)| +10|\Sigma _{(1,1)} (0)|=10, \notag \\
&b_{1} (X_{(1,1)}) = b_{3} (X_{(1,1)}) = b_{5} (X_{(1,1)}) = b_{7} (X_{(1,1)}) = b_{9} (X_{(1,1)}) =0.
\end{align}
Therefore, Poincar\'e polynomial of $X_{(1,1)}$ is given by,
\begin{equation}
P_{(1,1)} (t) = 1+5t^2 +10t^4 +10t^6 +5t^8 +t^{10} =(1+t^2 )^5.
\end{equation}

\item[(ii) \ ($(d_1 ,d_2 )=(2,1)$)]
Since $|\Sigma _{(2,1)} (0)|=1, |\Sigma _{(2,1)} (1)|=14, |\Sigma _{(2,1)} (2)|=84, |\Sigma _{(2,1)} (3)|=280$, Betti numbers are computed as follows:
\begin{align}
&b_{0} (X_{(2,1)}) = b_{14} (X_{(2,1)}) =|\Sigma _{(2,1)} (0)|=1, \notag \\
&b_{2} (X_{(2,1)}) = b_{12} (X_{(2,1)}) =|\Sigma _{(2,1)} (1)|-7|\Sigma _{(2,1)} (0)|=7, \notag \\
&b_{4} (X_{(2,1)}) = b_{10} (X_{(2,1)}) =|\Sigma _{(2,1)} (2)|-6|\Sigma _{(2,1)} (1)| +21|\Sigma _{(2,1)} (0)|=21, \notag \\
&b_{6} (X_{(2,1)}) = b_{8} (X_{(2,1)}) =|\Sigma _{(2,1)} (3)|-5|\Sigma _{(2,1)} (2)| +15|\Sigma _{(2,1)} (1)| -35|\Sigma _{(2,1)} (0)|=35, \notag \\
&b_{1} (X_{(2,1)}) = b_{3} (X_{(2,1)}) = \dots = b_{11} (X_{(2,1)}) = b_{13} (X_{(2,1)}) = 0.
\end{align}
Therefore, Poincar\'e polynomial of $X_{(2,1)}$ is
\begin{equation}
P_{(2,1)} (t) = 1+7t^2 +21t^4 +35t^6 +35t^8 + 21t^{10} +7t^{12} +t^{14} =(1+t^2 )^7.
\end{equation}
\end{description}
From these results, we are tempted to expect that $P_{(d_1 ,1)}(t)=(1+t^2)^{2d_1 +3}$, but it is not true. Since $|\Sigma _{(d_1 ,1)} (0)|=1, |\Sigma _{(d_1 ,1)} (1)|=4d_1 +6, |\Sigma _{(d_1 ,1)} (2)|=(15d_1 ^2 +43d_1 +22)/2, \dots$, some of Betti numbers of $X_{(d_1 ,1)}$ are
given by,
\begin{align}
&b_{0} (X_{(d_1 ,1)})=b_{4d_1 +6} (X_{(d_1 ,1)})=1, \notag \\
&b_{2} (X_{(d_1 ,1)})=b_{4d_1 +4} (X_{(d_1 ,1)})=2d_1 +3, \notag \\
&b_{4} (X_{(d_1 ,1)})=b_{4d_1 +2} (X_{(d_1 ,1)})=\frac{3d_1 ^2 +13d_1 +4}{2}, \notag \\
& \dots
\end{align}
Therefore, Poincar\'e polynomial of $X_{(d_1,1)}$ takes the form:
\begin{align}
P_{(d_1 ,1)}(t) =&1+(2d_1 +3)t^2 +\frac{3d_1 ^2 +13d_1 +4}{2} t^4 \notag \\ 
                   &+ \dots + \frac{3d_1 ^2 +13d_1 +4}{2} t^{4d_1 +2} + (2d_1 +3)t^{4d_1 +4} +t^{4d_1 +6}.
\end{align}
However, it is not equal to, 
\begin{align}
(1+t^2)^{2d_1 +3}=&1+(2d_1 +3)t^2 + (d_1 +1)(2d_2 +3)t^4 \notag \\ 
                       &+ \dots +(d_1 +1)(2d_2 +3)t^{4d_1 +2} +(2d_1 +3)t^{4d_1 +4} +t^{4d_1 +6}.
\end{align}
As we can see from the above, Poincar\'e polynomial in general $(d_{1},1)$ case might be more complicated.

\end{document}